# The Kumaraswamy Generalized Marshall-Olkin-G family of distributions

Laba Handique and Subrata Chakraborty*

Department of Statistics, Dibrugarh University

Dibrugarh-786004, India

*Corresponding Author. Email: subrata_stats@dibru.ac.in

**Abstract**

Another new family of continuous probability distribution is proposed by using Generalized Marshal-Olkin distribution as the base line distribution in the Kumaraswamy-*G* distribution. This family includes $Kw-G$ (Cordeiro and de Castro, 2011) and *GMO* (Jayakumar and Mathew, 2008) families special case besides a under of other distributions. The probability density function (pdf) and the survival function (sf) are expressed as series to observe as a mixture of the Generalized Marshal-Olkin distribution. Series expansions pdf of order statistics are also obtained. Moments, moment generating function, Rényi entropies, quantile function, random sample generation and asymptotes are also investigated. Parameter estimation by method of maximum likelihood and method of moment are also presented. Finally the proposed model is compared to the Generalized Marshall-Olkin Kumaraswamy extended family (Handique and Chakraborty, 2015) by considering four examples of real life data modeling.

## 1. Introduction

Recently, some efforts have been made to define new families of continuous distributions to extend well-known distributions and at the same time provide great flexibility in modelling data in practice. So, several classes by adding one or more parameters to generate new distributions have been proposed in the statistical literature. The basic motivation of these works is to bring in more flexibility in the modelling different type of data generated from real life situation.

Newly, there is renewed activity in this area to propose and investigate new families of distributions. A recent review paper by Tahir *et al.* (2015) provides a detail account of various popular techniques of generating new families of univariate continuous distributions through introduction of additional parameters. Very recent contributions in this line include the Kumaraswamy Marshal - Olkin family proposed by Alizadeh *et al.*, (2015), Marshal - Olkin Kumaraswamy-*G* family introduced by Handique and Chakraborty (2015a), Generalized Marshal - Olkin Kumaraswamy-*G* family introduced by Handique and Chakraborty (2015b) and Beta Generated Kumaraswamy-*G* family introduced by Handique and Chakraborty (2016a) among others.

In this article we propose another family of continuous probability distribution by integrating the Generalized Marshall Olkin family (Jayakumar and Mathew, 2008) as the base line distribution in the Kumaraswamy-*G* family (Cordeiro and de Castro, 2011). This new family referred to as the Kumaraswamy Generalized Marshall-Olkin-*G* ( $KwGMO-G$ ) family of distribution is investigated. We have discussed its mathematical properties, and applied it to a real data for illustration of the flexibility of the model.

The rest of this article is organized in seven sections. In section 2 we briefly introduce some important characteristics of probability distributions, the Kumaraswamy-*G* (Cordeiro and de Castro, 2011) and the generalized Marshall-Olkin family (Jayakumar and Mathew, 2008) of distributions. The proposed new family is defined along with its physical basis in section 3. Important special cases of the family along with their shape and main reliability characteristics are presented in the next section. In section 5 we discuss some general results of the proposed family, while different methods of estimation of parameters along with four comparative data modelling applications are presented in section 6. The article ends with a conclusion in section 7 followed by an appendix to derive asymptotic confidence bounds.

**2. Some formulas and notations**

Here first we list some formulas to be used in the subsequent sections of this article.

If $T$ is a continuous random variable with pdf, $f(t)$ and cdf $F(t)=P[T\le t]$, then its

Survival function (sf): $\bar{F}(t)=P[T>t]=1-F(t)$,

Hazard rate function (hrf): $h(t)=f(t)/\bar{F}(t)$,

Reverse hazard rate function (rhrf): $r(t)=f(t)/F(t)$,

Cumulative hazard rate function (chrf): $H(t)=-\log[\bar{F}(t)]$,

$(p,q,r)^{\text{th}}$ Power Weighted Moment (PWM): $\Gamma_{p,q,r}=\int_{-\infty}^{\infty}t^{p}[F(t)]^{q}[1-F(t)]^{r}f(t)\,dt$,

Rényi entropy: $I_R(\delta) = (1-\delta)^{-1} \log\left( \int_{-\infty}^{\infty} f(t)^{\delta}\, dt \right)$.

### 2.1 Kumaraswamy-*G* ( $Kw-G$ ) family of distribution

For a baseline cdf $F(t)$ with pdf $f(t)$, Cordeiro and de Castro (2011) defined $Kw-G$ distribution with cdf and pdf

$$F^{KwG}(t) = 1 - [1 - F(t)^a]^b, \quad 0 < t < \infty,\ 0 < a, b < \infty \tag{1}$$

and

$$f^{KwG}(t) = ab\, f(t) F(t)^{a-1} [1 - F(t)^a]^{b-1} \tag{2}$$

Where $t > 0,\ f(t) = \frac{d}{dt} F(t)$ and $a > 0, b > 0$ are shape parameters in addition to those in the baseline distribution. The sf, hrf, rhrf and chrf of this distribution are respectively given by

$$\overline{F}^{KwG}(t) = 1 - F^{KwG}(t) = [1 - F(t)^a]^b$$

$$h^{KwG}(t) = f^{KwG}(t) \big/ \overline{F}^{KwG}(t) = ab\, f(t) F(t)^{a-1} [1 - F(t)^a]^{b-1} \big/ [1 - F(t)^a]^b$$

$$= ab\, f(t) F(t)^{a-1} [1 - F(t)^a]^{-1}$$

$$r^{KwG}(t) = f^{KwG}(t) \big/ \overline{F}^{KwG}(t) = ab\, f(t) F(t)^{a-1} [1 - F(t)^a]^{b-1} \big/ 1 - [1 - F(t)^a]^b$$

$$= ab\, f(t) F(t)^{a-1} [1 - F(t)^a]^{b-1} \left\{ 1 - [1 - F(t)^a]^b \right\}^{-1}$$

and $H^{KwG}(t) = -b \log[1 - F(t)^a]$.

### 2.2 Generalized Marshall-Olkin Extended ( GMOE ) family of distribution

Jayakumar and Mathew (2008) proposed a generalization of the Marshall and Olkin (1997) family of distributions by using the Lehman second alternative (Lehmann 1953) to obtain the sf $\overline{F}^{GMO}(t)$ of the GMOE family of distributions by exponentiation the sf of MOE family of distributions as

$$\overline{F}^{GMO}(t) = \left[ \frac{\alpha \overline{G}(t)}{1 - \overline{\alpha}\, \overline{G}(t)} \right]^{\theta}, -\infty < t < \infty; 0 < \alpha < \infty; 0 < \theta < \infty \tag{3}$$

where $-\infty < t < \infty,\ \alpha > 0\ (\overline{\alpha} = 1 - \alpha)$ and $\theta > 0$ is an additional shape parameter. When $\theta = 1$, $\overline{F}^{GMO}(t) = \overline{F}^{MO}(t)$ and for $\alpha = \theta = 1$, $\overline{F}^{GMO}(t) = \overline{F}(t)$. The cdf and pdf of the GMOE distribution are respectively

$$F^{GMO}(t) = 1 - \left[ \frac{\alpha \overline{G}(t)}{1 - \overline{\alpha}\, \overline{G}(t)} \right]^{\theta} \tag{4}$$

and

$$f^{GMO}(t) = \theta \left[ \frac{\alpha \overline{G}(t)}{1 - \overline{\alpha}\, \overline{G}(t)} \right]^{\theta - 1} \left\{ \frac{\alpha\, g(t)}{[1 - \overline{\alpha}\, \overline{G}(t)]^2} \right\} = \frac{\theta \alpha^{\theta} g(t) \overline{G}(t)^{\theta - 1}}{[1 - \overline{\alpha}\, \overline{G}(t)]^{\theta + 1}} \tag{5}$$

Reliability measures like the hrf, rhrf and chrf associated with (1) are

$$h^{GMO}(t)=\frac{f^{GMO}(t)}{\overline{F}^{GMO}(t)}=\frac{\theta\alpha^{\theta}g(t)\overline{G}(t)^{\theta-1}}{[1-\overline{\alpha}\,\overline{G}(t)]^{\theta+1}}\Bigg/\left[\frac{\alpha\overline{G}(t)}{1-\overline{\alpha}\,\overline{G}(t)}\right]^{\theta}$$

$$=\theta\frac{g(t)}{\overline{G}(t)}\frac{1}{1-\overline{\alpha}\,\overline{G}(t)}=\theta\,\frac{h(t)}{1-\overline{\alpha}\,\overline{G}(t)}$$

$$r^{GMO}(t)=\frac{f^{GMO}(t)}{F^{GMO}(t)}=\frac{\theta\alpha^{\theta}g(t)\overline{G}(t)^{\theta-1}}{[1-\overline{\alpha}\,\overline{G}(t)]^{\theta+1}}\Bigg/1-\left[\frac{\alpha\overline{G}(t)}{1-\overline{\alpha}\,\overline{G}(t)}\right]^{\theta}$$

$$=\frac{\theta\alpha^{\theta}g(t)\overline{G}(t)^{\theta-1}}{[1-\overline{\alpha}\,\overline{G}(t)]\,[\{1-\overline{\alpha}\,\overline{G}(t)\}^{\theta}-\alpha^{\theta}\overline{G}(t)^{\theta}]}=\frac{\theta\alpha^{\theta}g(t)\overline{G}(t)^{\theta-1}}{[1-\overline{\alpha}\,\overline{G}(t)]^{\theta+1}-\alpha^{\theta}\overline{G}(t)^{\theta}[1-\overline{\alpha}\,\overline{G}(t)]}$$

$$H^{GMO}(t)=-\log\left[\frac{\alpha\overline{G}(t)}{1-\overline{\alpha}\,\overline{G}(t)}\right]^{\theta}=-\theta\log\left[\frac{\alpha\overline{G}(t)}{1-\overline{\alpha}\,\overline{G}(t)}\right]$$

Where $g(t)$, $G(t)$, $\overline{G}(t)$ and $h(t)$ are respectively the pdf, cdf, sf and hrf of the baseline distribution.

We denote the family of distribution with pdf (5) as $\text{GMOE}(\alpha,\theta,a,b)$ which for $\theta=1$, reduces to $\text{MOE}(\alpha,a,b)$.

**3. Kumaraswamy Generalized Marshall-Olkin-*G* ( $KwGMO-G$ ) family of distribution**

We now propose a new extension of the $Kw-G$ (Cordeiro and de Castro, 2011) family by considering the cdf and pdf of $GMO$ (Jayakumar and Mathew, 2008) distribution in (4) and (5) as the $f(t)$ and $F(t)$ respectively in the $Kw-G$ formulation in (2) and call it $KwGMO-G$ distribution. The resulting expression for the pdf of $KwGMO-G$ is given by

$$f^{KwGMO}(t;a,b,\alpha,\theta)$$

$$=\frac{ab\theta\alpha^{\theta}g(t)\overline{G}(t)^{\theta-1}}{[1-\overline{\alpha}\,\overline{G}(t)]^{\theta+1}}\left[1-\left\{\frac{\alpha\overline{G}(t)}{1-\overline{\alpha}\,\overline{G}(t)}\right\}^{\theta}\right]^{a-1}\left[1-\left[1-\left\{\frac{\alpha\overline{G}(t)}{1-\overline{\alpha}\,\overline{G}(t)}\right\}^{\theta}\right]^{a}\right]^{b-1} \quad (6)$$

$$,-\infty<t<\infty;0<\alpha<\infty;0<\theta<\infty,a>0,b>0$$

The cdf, sf, hrf, rhrf and chrf of $KwGMO-G$ distribution are respectively given by

cdf: $$F^{KwGMO}(t;a,b,\alpha,\theta)=1-\left[1-\left[1-\left\{\frac{\alpha\overline{G}(t)}{1-\overline{\alpha}\,\overline{G}(t)}\right\}^{\theta}\right]^{a}\right]^{b} \quad (7)$$

sf: $$\overline{F}^{KwGMO}(t;a,b,\alpha,\theta)=\left[1-\left[1-\left\{\frac{\alpha\overline{G}(t)}{1-\overline{\alpha}\,\overline{G}(t)}\right\}^{\theta}\right]^{a}\right]^{b}$$

hrf:

$$h^{KwGMO}(t;a,b,\alpha,\theta)=\frac{ab\theta\alpha^{\theta}\,g(t)\overline{G}(t)^{\theta-1}}{[1-\overline{\alpha}\,\overline{G}(t)]^{\theta+1}}\left[1-\left\{\frac{\alpha\overline{G}(t)}{1-\overline{\alpha}\,\overline{G}(t)}\right\}^{\theta}\right]^{a-1}\left[1-\left[1-\left\{\frac{\alpha\overline{G}(t)}{1-\overline{\alpha}\,\overline{G}(t)}\right\}^{\theta}\right]^{a}\right]^{-1} \quad (8)$$

rhrf:

$$r^{KwGMO}(t;a,b,\alpha,\theta)=\frac{ab\theta\alpha^{\theta}\,g(t)\overline{G}(t)^{\theta-1}}{[1-\overline{\alpha}\,\overline{G}(t)]^{\theta+1}}\left[1-\left\{\frac{\alpha\overline{G}(t)}{1-\overline{\alpha}\,\overline{G}(t)}\right\}^{\theta}\right]^{a-1}$$

$$\times\left[1-\left[1-\left\{\frac{\alpha\overline{G}(t)}{1-\overline{\alpha}\,\overline{G}(t)}\right\}^{\theta}\right]^{a}\right]^{b-1}\left[1-\left[1-\left[1-\left\{\frac{\alpha\overline{G}(t)}{1-\overline{\alpha}\,\overline{G}(t)}\right\}^{\theta}\right]^{a}\right]^{b}\right]^{-1} \quad (9)$$

chrf: $H^{KwGMO}(t;a,b,\alpha,\theta)=-b\log\left[1-\left[1-\left\{\frac{\alpha\overline{G}(t)}{1-\overline{\alpha}\,\overline{G}(t)}\right\}^{\theta}\right]^{a}\right]$

Remark: $KwGMO-G$ reduces to some known families of distributions as:

(i) For $\theta=1$, $f^{KwGMO}(t;a,b,\alpha,\theta)=f^{KwMO}(t;a,b,\alpha)$ (Alizadeh *et al*., 2015), (ii) for $\alpha=1$, $f^{KwGMO}(t;a,b,\alpha,\theta)=f^{GKwG}(t;a,b,\theta)$, (iii) for $a=b=1$, $f^{KwGMO}(t;a,b,\alpha,\theta)=f^{GMO}(t;\alpha,\theta)$ (Jayakumar and Mathew, 2008) and for $\alpha=\theta=1$, $f^{KwGMO}(t;a,b,\alpha,\theta)=f^{KwG}(t;a,b)$, (Cordeiro and de Castro, 2011).

### 3.1 Shape of the density and hazard function

Here we have plotted the pdf and hrf of the $KwGMO-G$ for some choices of the parameters to study the variety of shapes assumed by the family.

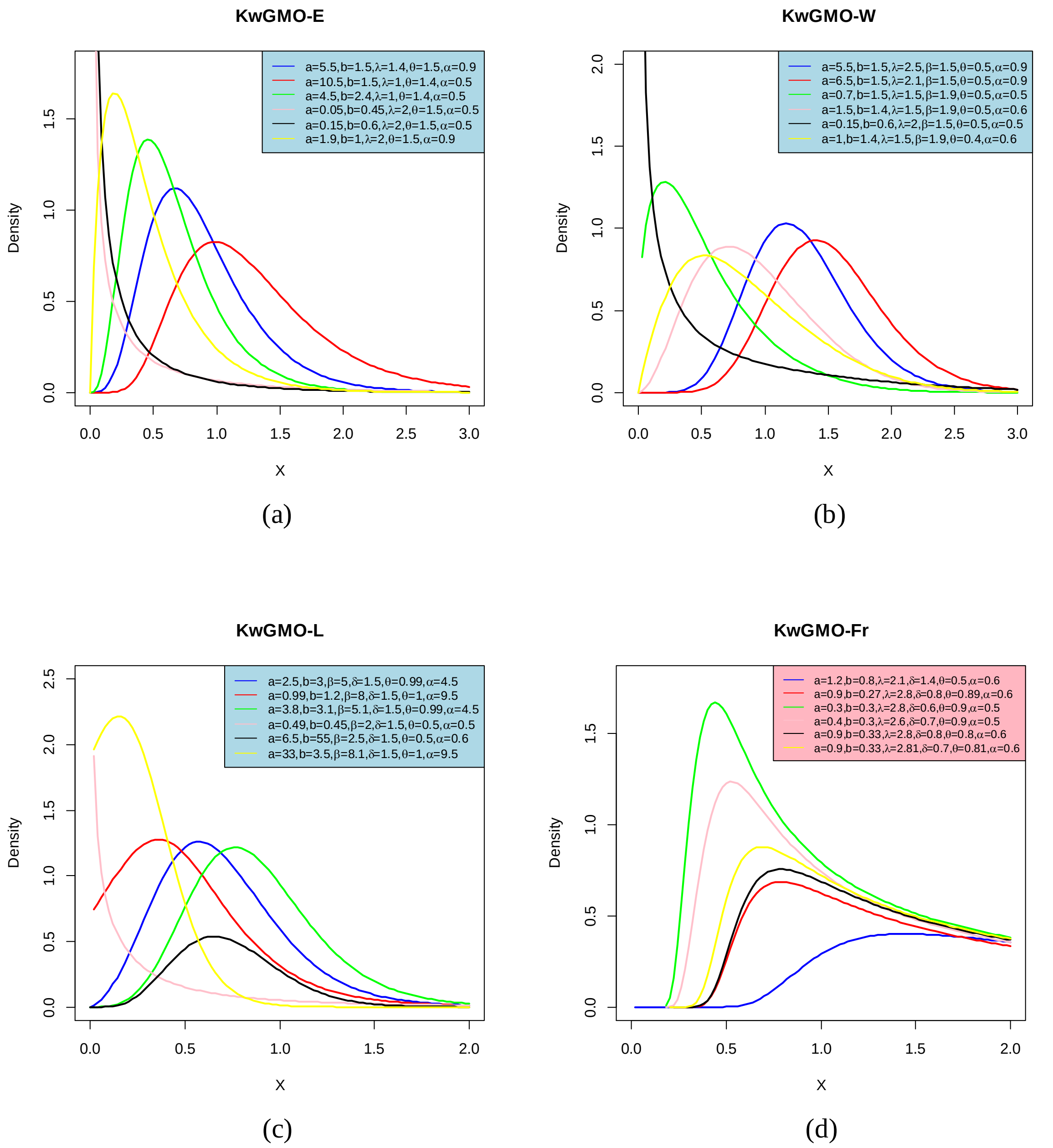

**Fig 1:** Density plots (a) $KwGMO-E$, (b) $KwGMO-W$, (c) $KwGMO-L$ and (d) $KwGMO-Fr$ distributions.

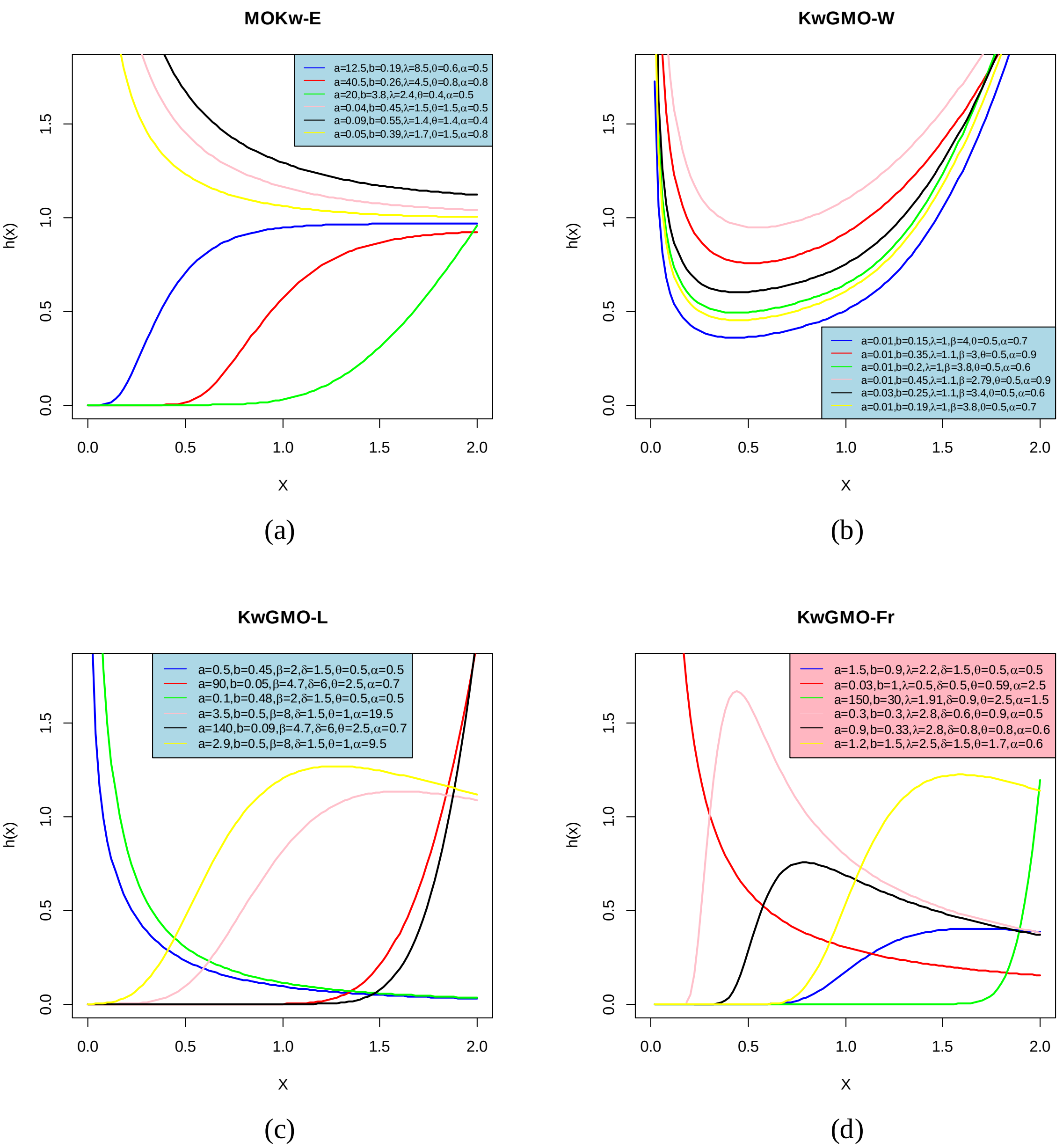

**Fig 2:** Hazard plots (a) $KwGMO-E$, (b) $KwGMO-W$, (c) $KwGMO-L$ and (d) $KwGMO-Fr$ distributions.

From the plots in figure 1 and 2 it can be seen that the family is very flexible and can offer many different types of shapes of density and hazard rate function including the bath tub shaped free hazard.

## 4. Some special $KwGMO-G$ distributions

In this section we provide some special cases of the $KwGMO-G$ family of distributions and list their main distributional characteristics.

4.1 The $KwGMO$ – exponential ( $KwGMO-E$ ) distribution:

Let the base line distribution be exponential with parameter $\lambda > 0$, $g(t:\lambda) = \lambda e^{-\lambda t}$, $t > 0$ and $G(t:\lambda) = 1 - e^{-\lambda t}$, $t > 0$ then for the $KwGMO-E$ model we get the pdf and cdf respectively as:

$$f^{KwGMOE}(t;a,b,\alpha,\theta,\lambda)$$

$$= \frac{ab\theta\alpha^{\theta}\lambda e^{-\lambda t}[e^{-\lambda t}]^{\theta-1}}{[1-\overline{\alpha}\,e^{-\lambda t}]^{\theta+1}}\left[1-\left\{\frac{\alpha e^{-\lambda t}}{1-\overline{\alpha}\,e^{-\lambda t}}\right\}^{\theta}\right]^{a-1}\left[1-\left[1-\left\{\frac{\alpha e^{-\lambda t}}{1-\overline{\alpha}\,e^{-\lambda t}}\right\}^{\theta}\right]^{a}\right]^{b-1}$$

cdf: $$F^{KwGMOE}(t;a,b,\alpha,\theta,\lambda) = 1-\left[1-\left[1-\left\{\frac{\alpha e^{-\lambda t}}{1-\overline{\alpha}\,e^{-\lambda t}}\right\}^{\theta}\right]^{a}\right]^{b}$$

sf: $$\overline{F}^{KwGMOE}(t;a,b,\alpha,\theta,\lambda) = \left[1-\left[1-\left\{\frac{\alpha e^{-\lambda t}}{1-\overline{\alpha}\,e^{-\lambda t}}\right\}^{\theta}\right]^{a}\right]^{b}$$

hrf: $h^{KwGMOE}(t;a,b,\alpha,\theta,\lambda)$

$$= \frac{ab\theta\alpha^{\theta}\lambda e^{-\lambda t}[e^{-\lambda t}]^{\theta-1}}{[1-\overline{\alpha}\,e^{-\lambda t}]^{\theta+1}}\left[1-\left\{\frac{\alpha e^{-\lambda t}}{1-\overline{\alpha}\,e^{-\lambda t}}\right\}^{\theta}\right]^{a-1}\left[1-\left[1-\left\{\frac{\alpha e^{-\lambda t}}{1-\overline{\alpha}\,e^{-\lambda t}}\right\}^{\theta}\right]^{a}\right]^{-1}$$

rhrf: $$r^{KwGMOE}(t;a,b,\alpha,\theta) = \frac{ab\theta\alpha^{\theta}\lambda e^{-\lambda t}(e^{-\lambda t})^{\theta-1}}{[1-\overline{\alpha}\,e^{-\lambda t}]^{\theta+1}}\left[1-\left\{\frac{\alpha e^{-\lambda t}}{1-\overline{\alpha}\,e^{-\lambda t}}\right\}^{\theta}\right]^{a-1}$$

$$\left[1-\left[1-\left\{\frac{\alpha e^{-\lambda t}}{1-\overline{\alpha}\,e^{-\lambda t}}\right\}^{\theta}\right]^{a}\right]^{b-1}\left[1-\left[1-\left[1-\left\{\frac{\alpha e^{-\lambda t})}{1-\overline{\alpha}\,e^{-\lambda t})}\right\}^{\theta}\right]^{a}\right]^{b}\right]^{-1}$$

chrf: $$H^{KwGMOE}(t;a,b,\alpha,\theta) = -b\log\left[1-\left[1-\left\{\frac{\alpha e^{-\lambda t}}{1-\overline{\alpha}\,e^{-\lambda t}}\right\}^{\theta}\right]^{a}\right]$$

### 4.2 The $KwGMO-$ Lomax ( $KwGMO-L$ ) distribution

Considering the Lomax distribution (Ghitany *et al.* 2007) with pdf and cdf given by $g(t:\beta,\delta) = (\beta/\delta)[1+(t/\delta)]^{-(\beta+1)}$, $t > 0$, and $G(t:\beta,\delta) = 1-[1+(t/\delta)]^{-\beta}$, $\beta > 0$ and $\delta > 0$ the pdf and cdf of the $KwGMO-L$ distribution are given by

$$f^{KwGMOL}(t;a,b,\alpha,\theta,\beta,\delta) = \frac{ab\theta\alpha^{\theta}(\beta/\delta)[1+(t/\delta)]^{-(\beta+1)}[\{1+(t/\delta)\}^{-\beta}]^{\theta-1}}{[1-\overline{\alpha}\{1+(t/\delta)\}^{-\beta}]^{\theta+1}}$$

$$\left[1-\left\{\frac{\alpha[1+(t/\delta)]^{-\beta}}{1-\bar{\alpha}[1+(t/\delta)]^{-\beta}}\right\}^{\theta}\right]^{a-1}\left[1-\left[1-\left\{\frac{\alpha[1+(t/\delta)]^{-\beta}}{1-\bar{\alpha}[1+(t/\delta)]^{-\beta}}\right\}^{\theta}\right]^{a}\right]^{b-1}$$

cdf: $$F^{KwGMOL}(t;a,b,\alpha,\theta,\beta,\delta)=1-\left[1-\left[1-\left\{\frac{\alpha[1+(t/\delta)]^{-\beta}}{1-\bar{\alpha}[1+(t/\delta)]^{-\beta}}\right\}^{\theta}\right]^{a}\right]^{b}$$

sf: $$\bar{F}^{KwGMOL}(t;a,b,\alpha,\theta,\beta,\delta)=\left[1-\left[1-\left\{\frac{\alpha[1+(t/\delta)]^{-\beta}}{1-\bar{\alpha}[1+(t/\delta)]^{-\beta}}\right\}^{\theta}\right]^{a}\right]^{b}$$

hrf:

$$h^{KwGMOL}(t;a,b,\alpha,\theta,\beta,\delta)=\frac{ab\theta\alpha^{\theta}(\beta/\delta)[1+(t/\delta)]^{-(\beta+1)}\bar{G}(t)^{\theta-1}}{[1-\bar{\alpha}[1+(t/\delta)]^{-\beta}]^{\theta+1}}\left[1-\left\{\frac{\alpha[1+(t/\delta)]^{-\beta}}{1-\bar{\alpha}[1+(t/\delta)]^{-\beta}}\right\}^{\theta}\right]^{a-1}$$

$$\left[1-\left[1-\left\{\frac{\alpha[1+(t/\delta)]^{-\beta}}{1-\bar{\alpha}[1+(t/\delta)]^{-\beta}}\right\}^{\theta}\right]^{a}\right]^{-1}$$

rhrf: $r^{KwGMOL}(t;a,b,\alpha,\theta,\beta,\delta)$

$$=\frac{ab\theta\alpha^{\theta}(\beta/\delta)[1+(t/\delta)]^{-(\beta+1)}[[1+(t/\delta)]^{-\beta}]^{\theta-1}}{[1-\bar{\alpha}[1+(t/\delta)]^{-\beta}]^{\theta+1}}\left[1-\left\{\frac{\alpha[1+(t/\delta)]^{-\beta}}{1-\bar{\alpha}[1+(t/\delta)]^{-\beta}}\right\}^{\theta}\right]^{a-1}$$

$$\left[1-\left[1-\left\{\frac{\alpha[1+(t/\delta)]^{-\beta}}{1-\bar{\alpha}[1+(t/\delta)]^{-\beta}}\right\}^{\theta}\right]^{a}\right]^{b-1}\left[1-\left[1-\left[1-\left\{\frac{\alpha[1+(t/\delta)]^{-\beta}}{1-\bar{\alpha}[1+(t/\delta)]^{-\beta}}\right\}^{\theta}\right]^{a}\right]^{b}\right]^{-1}$$

chrf: $$r^{KwGMOL}(t;a,b,\alpha,\theta,\beta,\delta)=-b\log\left[1-\left[1-\left\{\frac{\alpha[1+(t/\delta)]^{-\beta}}{1-\bar{\alpha}[1+(t/\delta)]^{-\beta}}\right\}^{\theta}\right]^{a}\right]$$

### 4.3 The $KwGMO-$Weibull ($KwGMO-W$) distribution

Considering the Weibull distribution (Ghitany *et al.* 2005, Zhang and Xie 2007) with parameters $\lambda>0$ and $\beta>0$ having pdf and cdf $g(t)=\lambda\beta\, t^{\beta-1}e^{-\lambda t^{\beta}}$ and $G(t)=1-e^{-\lambda t^{\beta}}$ respectively we get the pdf and cdf of $KwGMO-W$ distribution as

$$f^{KwGMOW}(t;a,b,\alpha,\theta,\beta,\lambda)=\frac{ab\theta\alpha^{\theta}\lambda\beta\, t^{\beta-1}e^{-\lambda t^{\beta}}[e^{-\lambda t^{\beta}}]^{\theta-1}}{[1-\bar{\alpha}\,e^{-\lambda t^{\beta}}]^{\theta+1}}\left[1-\left\{\frac{\alpha\, e^{-\lambda t^{\beta}}}{1-\bar{\alpha}\,e^{-\lambda t^{\beta}}}\right\}^{\theta}\right]^{a-1}$$

$$\left[1-\left[1-\left\{\frac{\alpha\, e^{-\lambda t^{\beta}}}{1-\bar{\alpha}\, e^{-\lambda t^{\beta}}}\right\}^{\theta}\right]^{a}\right]^{b-1}$$

cdf: $$F^{KwGMOW}(t;a,b,\alpha,\theta,\beta,\lambda)=1-\left[1-\left[1-\left\{\frac{\alpha\, e^{-\lambda t^{\beta}}}{1-\bar{\alpha}\, e^{-\lambda t^{\beta}}}\right\}^{\theta}\right]^{a}\right]^{b}$$

sf: $$\bar{F}^{KwGMOW}(t;a,b,\alpha,\theta,\beta,\lambda)=\left[1-\left[1-\left\{\frac{\alpha\, e^{-\lambda t^{\beta}}}{1-\bar{\alpha}\, e^{-\lambda t^{\beta}}}\right\}^{\theta}\right]^{a}\right]^{b}$$

hrf:

$$h^{KwGMOW}(t;a,b,\alpha,\theta,\beta,\lambda)=\frac{a\,b\,\theta\,\alpha^{\theta}\,\lambda\,\beta\, t^{\beta-1}\, e^{-\lambda t^{\beta}}\,[e^{-\lambda t^{\beta}}]^{\theta-1}}{[1-\bar{\alpha}\, e^{-\lambda t^{\beta}}]^{\theta+1}}\left[1-\left\{\frac{\alpha\, e^{-\lambda t^{\beta}}}{1-\bar{\alpha}\, e^{-\lambda t^{\beta}}}\right\}^{\theta}\right]^{a-1}$$

$$\left[1-\left[1-\left\{\frac{\alpha\, e^{-\lambda t^{\beta}}}{1-\bar{\alpha}\, e^{-\lambda t^{\beta}}}\right\}^{\theta}\right]^{a}\right]^{-1}$$

rhrf:

$$r^{KwGMOW}(t;a,b,\alpha,\theta,\beta,\lambda)=\frac{a\,b\,\theta\,\alpha^{\theta}\,\lambda\,\beta\, t^{\beta-1}\, e^{-\lambda t^{\beta}}\,[e^{-\lambda t^{\beta}}]^{\theta-1}}{[1-\bar{\alpha}\, e^{-\lambda t^{\beta}}]^{\theta+1}}\left[1-\left\{\frac{\alpha\, e^{-\lambda t^{\beta}}}{1-\bar{\alpha}\, e^{-\lambda t^{\beta}}}\right\}^{\theta}\right]^{a-1}$$

$$\left[1-\left[1-\left\{\frac{\alpha\, e^{-\lambda t^{\beta}}}{1-\bar{\alpha}\, e^{-\lambda t^{\beta}}}\right\}^{\theta}\right]^{a}\right]^{b-1}\left[1-\left[1-\left[1-\left\{\frac{\alpha\, e^{-\lambda t^{\beta}}}{1-\bar{\alpha}\, e^{-\lambda t^{\beta}}}\right\}^{\theta}\right]^{a}\right]^{b}\right]^{-1}$$

chrf: $$H^{KwGMOW}(t;a,b,\alpha,\theta,\beta,\lambda)=-b\log\left[1-\left[1-\left\{\frac{\alpha\, e^{-\lambda t^{\beta}}}{1-\bar{\alpha}\, e^{-\lambda t^{\beta}}}\right\}^{\theta}\right]^{a}\right]$$

### 4.4 The $KwGMO-$Frechet ( $KwGMO-Fr$ ) distribution

Suppose the base line distribution is the Frechet distribution (Krishna *et al*., 2013) with pdf and cdf given by $g(t)=\lambda\,\delta^{\lambda}\,t^{-(\lambda+1)}\,e^{-(\delta/t)^{\lambda}}$ and $G(t)=e^{-(\delta/t)^{\lambda}}$, $t>0$ respectively, and then the corresponding pdf and cdf of $KwGMO-Fr$ distribution becomes

$$f^{KwGMOFr}(t;a,b,\alpha,\theta,\lambda,\delta)=\frac{a\,b\,\theta\,\alpha^{\theta}\,\lambda\,\delta^{\lambda}\,t^{-(\lambda+1)}\,e^{-(\delta/t)^{\lambda}}\,[1-e^{-(\delta/t)^{\lambda}}]^{\theta-1}}{[1-\bar{\alpha}\,e^{-(\delta/t)^{\lambda}}]^{\theta+1}}$$

$$\left[1-\left\{\frac{\alpha[1-e^{-(\delta/t)^{\lambda}}]}{1-\bar{\alpha}[1-e^{-(\delta/t)^{\lambda}}]}\right\}^{\theta}\right]^{a-1}\left[1-\left[1-\left\{\frac{\alpha[1-e^{-(\delta/t)^{\lambda}}]}{1-\bar{\alpha}[1-e^{-(\delta/t)^{\lambda}}]}\right\}^{\theta}\right]^{a}\right]^{b-1}$$

cdf: $$F^{KwGMOFr}(t;a,b,\alpha,\theta,\lambda,\delta)=1-\left[1-\left[1-\left\{\frac{\alpha[1-e^{-(\delta/t)^{\lambda}}]}{1-\bar{\alpha}[1-e^{-(\delta/t)^{\lambda}}]}\right\}^{\theta}\right]^{a}\right]^{b}$$

sf: $$\bar{F}^{KwGMOFr}(t;a,b,\alpha,\theta,\lambda,\delta)=\left[1-\left[1-\left\{\frac{\alpha[1-e^{-(\delta/t)^{\lambda}}]}{1-\bar{\alpha}[1-e^{-(\delta/t)^{\lambda}}]}\right\}^{\theta}\right]^{a}\right]^{b}$$

hrf: $$h^{KwGMOFr}(t;a,b,\alpha,\theta,\lambda,\delta)=\frac{ab\theta\alpha^{\theta}\lambda\delta^{\lambda}t^{-(\lambda+1)}e^{-(\delta/t)^{\lambda}}[1-e^{-(\delta/t)^{\lambda}}]^{\theta-1}}{[1-\bar{\alpha}(1-e^{-(\delta/t)^{\lambda}})]^{\theta+1}}$$

$$\left[1-\left\{\frac{\alpha[1-e^{-(\delta/t)^{\lambda}}]}{1-\bar{\alpha}[1-e^{-(\delta/t)^{\lambda}}]}\right\}^{\theta}\right]^{a-1}\left[1-\left[1-\left\{\frac{\alpha[1-e^{-(\delta/t)^{\lambda}}]}{1-\bar{\alpha}[1-e^{-(\delta/t)^{\lambda}}]}\right\}^{\theta}\right]^{a}\right]^{-1}$$

rhrf: $r^{KwGMOFr}(t;a,b,\alpha,\theta,\lambda,\delta)=$

$$\frac{ab\theta\alpha^{\theta}\lambda\delta^{\lambda}t^{-(\lambda+1)}e^{-(\delta/t)^{\lambda}}[1-e^{-(\delta/t)^{\lambda}}]^{\theta-1}}{[1-\bar{\alpha}1-e^{-(\delta/t)^{\lambda}}]^{\theta+1}}\left[1-\left\{\frac{\alpha[1-e^{-(\delta/t)^{\lambda}}]}{1-\bar{\alpha}[1-e^{-(\delta/t)^{\lambda}}]}\right\}^{\theta}\right]^{a-1}$$

$$\left[1-\left[1-\left\{\frac{\alpha[1-e^{-(\delta/t)^{\lambda}}]}{1-\bar{\alpha}[1-e^{-(\delta/t)^{\lambda}}]}\right\}^{\theta}\right]^{a}\right]^{b-1}\left[1-\left[1-\left[1-\left\{\frac{\alpha[1-e^{-(\delta/t)^{\lambda}}]}{1-\bar{\alpha}[1-e^{-(\delta/t)^{\lambda}}]}\right\}^{\theta}\right]^{a}\right]^{b}\right]^{-1}$$

chrf: $$H^{KwGMOFr}(t;a,b,\alpha,\theta,\lambda,\delta)=-b\log\left[1-\left[1-\left\{\frac{\alpha[1-e^{-(\delta/t)^{\lambda}}]}{1-\bar{\alpha}[1-e^{-(\delta/t)^{\lambda}}]}\right\}^{\theta}\right]^{a}\right]^{b}$$

**4.5 The $KwGMO-$Gompertz ($KwGMO-Go$) distribution**

Next by taking the Gompertz distribution (Gieser et al. 1998) with pdf and cdf $g(t)=\beta e^{\lambda t}e^{-\frac{\beta}{\lambda}(e^{\lambda t}-1)}$ and $G(t)=1-e^{-\frac{\beta}{\lambda}(e^{\lambda t}-1)}$, $\beta>0,\ \lambda>0,\ t>0$ respectively, we get the pdf and cdf of $KwGMO-Go$ distribution as

$$f^{KwGMOGo}(t;a,b,\alpha,\theta,\beta,\lambda)=\frac{ab\theta\alpha^{\theta}\beta e^{\lambda t}e^{-\frac{\beta}{\lambda}(e^{\lambda t}-1)}[e^{-\frac{\beta}{\lambda}(e^{\lambda t}-1)}]^{\theta-1}}{[1-\bar{\alpha}e^{-\frac{\beta}{\lambda}(e^{\lambda t}-1)}]^{\theta+1}}\left[1-\left\{\frac{\alpha e^{-\frac{\beta}{\lambda}(e^{\lambda t}-1)}}{1-\bar{\alpha}e^{-\frac{\beta}{\lambda}(e^{\lambda t}-1)}}\right\}^{\theta}\right]^{a-1}$$

$$\times\left[1-\left[1-\left\{\frac{\alpha e^{-\frac{\beta}{\lambda}(e^{\lambda t}-1)}}{1-\bar{\alpha} e^{-\frac{\beta}{\lambda}(e^{\lambda t}-1)}}\right\}^{\theta}\right]^{a}\right]^{b-1}$$

cdf: $$F^{KwGMOGo}(t;a,b,\alpha,\theta,,\beta,\lambda)=1-\left[1-\left[1-\left\{\frac{\alpha e^{-\frac{\beta}{\lambda}(e^{\lambda t}-1)}}{1-\bar{\alpha} e^{-\frac{\beta}{\lambda}(e^{\lambda t}-1)}}\right\}^{\theta}\right]^{a}\right]^{b}$$

sf: $$\bar{F}^{KwGMOGo}(t;a,b,\alpha,\theta,,\beta,\lambda)=\left[1-\left[1-\left\{\frac{\alpha e^{-\frac{\beta}{\lambda}(e^{\lambda t}-1)}}{1-\bar{\alpha} e^{-\frac{\beta}{\lambda}(e^{\lambda t}-1)}}\right\}^{\theta}\right]^{a}\right]^{b}$$

hrf: $$h^{KwGMOGo}(t;a,b,\alpha,\theta,,\beta,\lambda)=\frac{ab\theta\alpha^{\theta}\beta e^{\lambda t} e^{-\frac{\beta}{\lambda}(e^{\lambda t}-1)}[e^{-\frac{\beta}{\lambda}(e^{\lambda t}-1)}]^{\theta-1}}{[1-\bar{\alpha} e^{-\frac{\beta}{\lambda}(e^{\lambda t}-1)}]^{\theta+1}}$$

$$\times\left[1-\left\{\frac{\alpha e^{-\frac{\beta}{\lambda}(e^{\lambda t}-1)}}{1-\bar{\alpha} e^{-\frac{\beta}{\lambda}(e^{\lambda t}-1)}}\right\}^{\theta}\right]^{a-1}\left[1-\left[1-\left\{\frac{\alpha e^{-\frac{\beta}{\lambda}(e^{\lambda t}-1)}}{1-\bar{\alpha} e^{-\frac{\beta}{\lambda}(e^{\lambda t}-1)}}\right\}^{\theta}\right]^{a}\right]^{-1}$$

rhrf:

$$r^{KwGMOGo}(t;a,b,\alpha,\theta,,\beta,\lambda)=\frac{ab\theta\alpha^{\theta}\beta e^{\lambda t} e^{-\frac{\beta}{\lambda}(e^{\lambda t}-1)}[e^{-\frac{\beta}{\lambda}(e^{\lambda t}-1)}]^{\theta-1}}{[1-\bar{\alpha} e^{-\frac{\beta}{\lambda}(e^{\lambda t}-1)}]^{\theta+1}}\left[1-\left\{\frac{\alpha e^{-\frac{\beta}{\lambda}(e^{\lambda t}-1)}}{1-\bar{\alpha} e^{-\frac{\beta}{\lambda}(e^{\lambda t}-1)}}\right\}^{\theta}\right]^{a-1}$$

$$\times\left[1-\left[1-\left\{\frac{\alpha e^{-\frac{\beta}{\lambda}(e^{\lambda t}-1)}}{1-\bar{\alpha} e^{-\frac{\beta}{\lambda}(e^{\lambda t}-1)}}\right\}^{\theta}\right]^{a}\right]^{b-1}\left[1-\left[1-\left[1-\left\{\frac{\alpha e^{-\frac{\beta}{\lambda}(e^{\lambda t}-1)}}{1-\bar{\alpha} e^{-\frac{\beta}{\lambda}(e^{\lambda t}-1)}}\right\}^{\theta}\right]^{a}\right]^{b}\right]^{-1}$$

chrf: $$H^{KwGMOGo}(t;a,b,\alpha,\theta,,\beta,\lambda)=-b\log\left[1-\left[1-\left\{\frac{\alpha e^{-\frac{\beta}{\lambda}(e^{\lambda t}-1)}}{1-\bar{\alpha} e^{-\frac{\beta}{\lambda}(e^{\lambda t}-1)}}\right\}^{\theta}\right]^{a}\right]$$

### 4.6 The $KwGMO-$Extended Weibull ($KwGMO-EW$) distribution

The pdf and the cdf of the extended Weibull (*EW*) distributions of Gurvich *et al.* (1997) is given by $g(t:\delta,\vartheta)=\delta\exp[-\delta Z(t:\vartheta)]z(t:\vartheta)$ and $G(t:\delta,\vartheta)=1-\exp[-\delta Z(t:\vartheta)]$, $t \in D \subseteq R_+, \delta > 0$ where $Z(t:\vartheta)$ is a non-negative monotonically increasing function which depends on the parameter vector $\vartheta$. and $z(t:\vartheta)$ is the derivative of $Z(t:\vartheta)$.

By considering $EW$ as the base line distribution we derive pdf and cdf of the $KwGMO-EW$ as

$$f^{KwGMOEW}(t;a,b,\alpha,\theta,\delta,\vartheta) = \frac{ab\theta\alpha^{\theta}\delta\exp[-\delta Z(t:\vartheta)]z(t:\vartheta)[\exp[-\delta Z(t:\vartheta)]]^{\theta-1}}{[1-\bar{\alpha}\exp[-\delta Z(t:\vartheta)]]^{\theta+1}}$$

$$\left[1-\left\{\frac{\alpha\exp[-\delta Z(t:\vartheta)]}{1-\bar{\alpha}\exp[-\delta Z(t:\vartheta)]}\right\}^{\theta}\right]^{a-1}\left[1-\left[1-\left\{\frac{\alpha\exp[-\delta Z(t:\vartheta)]}{1-\bar{\alpha}\exp[-\delta Z(t:\vartheta)]}\right\}^{\theta}\right]^{a}\right]^{b-1}$$

cdf: $$F^{KwGMOEW}(t;a,b,\alpha,\theta,\delta,\vartheta) = 1-\left[1-\left[1-\left\{\frac{\alpha\exp[-\delta Z(t:\vartheta)]}{1-\bar{\alpha}\exp[-\delta Z(t:\vartheta)]}\right\}^{\theta}\right]^{a}\right]^{b}$$

sf: $$\bar{F}^{KwGMOEW}(t;a,b,\alpha,\theta,\delta,\vartheta) = \left[1-\left[1-\left\{\frac{\alpha\exp[-\delta Z(t:\vartheta)]}{1-\bar{\alpha}\exp[-\delta Z(t:\vartheta)]}\right\}^{\theta}\right]^{a}\right]^{b}$$

hrf:

$$h^{KwGMOEW}(t;a,b,\alpha,\theta,\delta,\vartheta) = \frac{ab\theta\alpha^{\theta}\delta\exp[-\delta Z(t:\vartheta)]z(t:\vartheta)[\exp[-\delta Z(t:\vartheta)]]^{\theta-1}}{[1-\bar{\alpha}\exp[-\delta Z(t:\vartheta)]]^{\theta+1}}$$

$$\left[1-\left\{\frac{\alpha\exp[-\delta Z(t:\vartheta)]}{1-\bar{\alpha}\exp[-\delta Z(t:\vartheta)]}\right\}^{\theta}\right]^{a-1}\left[1-\left[1-\left\{\frac{\alpha\exp[-\delta Z(t:\vartheta)]}{1-\bar{\alpha}\exp[-\delta Z(t:\vartheta)]}\right\}^{\theta}\right]^{a}\right]^{-1}$$

rhr: $$r^{KwGMOEW}(t;a,b,\alpha,\theta,\delta,\vartheta) = \frac{ab\theta\alpha^{\theta}\delta\exp[-\delta Z(t:\vartheta)]z(t:\vartheta)[\exp[-\delta Z(t:\vartheta)]]^{\theta-1}}{[1-\bar{\alpha}\exp[-\delta Z(t:\vartheta)]]^{\theta+1}}$$

$$\left[1-\left\{\frac{\alpha\exp[-\delta Z(t:\vartheta)]}{1-\bar{\alpha}\exp[-\delta Z(t:\vartheta)]}\right\}^{\theta}\right]^{a-1}\left[1-\left[1-\left\{\frac{\alpha\exp[-\delta Z(t:\vartheta)]}{1-\bar{\alpha}\exp[-\delta Z(t:\vartheta)]}\right\}^{\theta}\right]^{a}\right]^{b-1}$$

$$\left[1-\left[1-\left[1-\left\{\frac{\alpha\exp[-\delta Z(t:\vartheta)]}{1-\bar{\alpha}\exp[-\delta Z(t:\vartheta)]}\right\}^{\theta}\right]^{a}\right]^{b}\right]^{-1}$$

chrf: $$H^{KwGMOEW}(t;a,b,\alpha,\theta,\delta,\vartheta) = -\log\left[1-\left[1-\left\{\frac{\alpha\exp[-\delta Z(t:\vartheta)]}{1-\bar{\alpha}\exp[-\delta Z(t:\vartheta)]}\right\}^{\theta}\right]^{a}\right]$$

Important models can be seen as particular cases with different choices of $Z(t:\vartheta)$:

(1) $Z(t:\vartheta) = t$: exponential distribution.

(2) $Z(t:\vartheta)=t^2$: Rayleigh (Burr type-X) distribution.

(3) $Z(t:\vartheta)=\log(t/k)$: Pareto distribution

(4) $Z(t:\vartheta)=\beta^{-1}[\exp(\beta t)-1]$: Gompertz distribution.

**4.7 The $KwGMO-$Extended Modified Weibull ($KwGMO-EMW$) distribution**

The modified Weibull (*MW*) distribution (Sarhan and Zaindin 2013) with cdf and pdf is given by

$G(t;\sigma,\beta,\gamma)=1-\exp[-\sigma t-\beta t^{\gamma}]$, $t>0, \gamma>0, \sigma,\beta\geq 0, \sigma+\beta>0$ and

$g(t;\sigma,\beta,\gamma)=(\sigma+\beta\gamma t^{\gamma-1})\exp[-\sigma t-\beta t^{\gamma}]$ respectively.

The corresponding pdf and cdf of $KwGMO-EMW$ are given by

$$f^{KwGMOEMW}(t;a,b,\alpha,\theta,\sigma,\beta,\gamma)=\frac{a b\theta\alpha^{\theta}(\sigma+\beta\gamma t^{\gamma-1})\exp[-\sigma t-\beta t^{\gamma}][\exp[-\sigma t-\beta t^{\gamma}]^{\theta-1}}{[1-\overline{\alpha}[\exp[-\sigma t-\beta t^{\gamma}]]^{\theta+1}}$$

$$\left[1-\left\{\frac{\alpha[\exp[-\sigma t-\beta t^{\gamma}]}{1-\overline{\alpha}[\exp[-\sigma t-\beta t^{\gamma}]}\right\}^{\theta}\right]^{a-1}\left[1-\left[1-\left\{\frac{\alpha[\exp[-\sigma t-\beta t^{\gamma}]}{1-\overline{\alpha}[\exp[-\sigma t-\beta t^{\gamma}]}\right\}^{\theta}\right]^{a}\right]^{b-1}$$

The cdf, sf, hrf, rhrf and chrf of $KwGMO-EMW$ distribution are respectively given by

cdf: $$F^{KwGMOEMW}(t;a,b,\alpha,\theta,\sigma,\beta,\gamma)=1-\left[1-\left[1-\left\{\frac{\alpha[\exp[-\sigma t-\beta t^{\gamma}]}{1-\overline{\alpha}[\exp[-\sigma t-\beta t^{\gamma}]}\right\}^{\theta}\right]^{a}\right]^{b}$$

sf: $$\overline{F}^{KwGMOEMW}(t;a,b,\alpha,\theta,\sigma,\beta,\gamma)=\left[1-\left[1-\left\{\frac{\alpha[\exp[-\sigma t-\beta t^{\gamma}]}{1-\overline{\alpha}[\exp[-\sigma t-\beta t^{\gamma}]}\right\}^{\theta}\right]^{a}\right]^{b}$$

hrf:

$$h^{KwGMOEMW}(t;a,b,\alpha,\theta,\sigma,\beta,\gamma)=\frac{a b\theta\alpha^{\theta}(\sigma+\beta\gamma t^{\gamma-1})\exp[-\sigma t-\beta t^{\gamma}][\exp[-\sigma t-\beta t^{\gamma}]^{\theta-1}}{[1-\overline{\alpha}[\exp[-\sigma t-\beta t^{\gamma}]]^{\theta+1}}$$

$$\left[1-\left\{\frac{\alpha[\exp[-\sigma t-\beta t^{\gamma}]}{1-\overline{\alpha}[\exp[-\sigma t-\beta t^{\gamma}]}\right\}^{\theta}\right]^{a-1}\left[1-\left[1-\left\{\frac{\alpha[\exp[-\sigma t-\beta t^{\gamma}]}{1-\overline{\alpha}[\exp[-\sigma t-\beta t^{\gamma}]}\right\}^{\theta}\right]^{a}\right]^{-1}$$

rhrf:

$$r^{KwGMOEMW}(t;a,b,\alpha,\theta,\sigma,\beta,\gamma)=\frac{a b\theta\alpha^{\theta}(\sigma+\beta\gamma t^{\gamma-1})\exp[-\sigma t-\beta t^{\gamma}][\exp[-\sigma t-\beta t^{\gamma}]^{\theta-1}}{[1-\overline{\alpha}[\exp[-\sigma t-\beta t^{\gamma}]]^{\theta+1}}$$

$$\left[1-\left\{\frac{\alpha[\exp[-\sigma t-\beta t^{\gamma}]}{1-\bar{\alpha}[\exp[-\sigma t-\beta t^{\gamma}]}\right\}^{\theta}\right]^{a-1}\left[1-\left[1-\left\{\frac{\alpha[\exp[-\sigma t-\beta t^{\gamma}]}{1-\bar{\alpha}[\exp[-\sigma t-\beta t^{\gamma}]}\right\}^{\theta}\right]^{a}\right]^{b-1}$$

$$\left[1-\left[1-\left[1-\left\{\frac{\alpha[\exp[-\sigma t-\beta t^{\gamma}]}{1-\bar{\alpha}[\exp[-\sigma t-\beta t^{\gamma}]}\right\}^{\theta}\right]^{a}\right]^{b}\right]^{-1}$$

chrf: $H^{KwGMOEMW}(t;a,b,\alpha,\theta,\sigma,\beta,\gamma)=-b\log\left[1-\left[1-\left\{\frac{\alpha[\exp[-\sigma t-\beta t^{\gamma}]}{1-\bar{\alpha}[\exp[-\sigma t-\beta t^{\gamma}]}\right\}^{\theta}\right]^{a}\right]$

### 4.8 The $KwGMO-$ Extended Exponentiated Pareto ( $KwGMO-EEP$ ) distribution

The pdf and cdf of the exponentiated Pareto distribution, of Nadarajah (2005), are given respectively by $g(t)=\gamma k\beta^{k}t^{-(k+1)}[1-(\beta/t)^{k}]^{\gamma-1}$ and $G(t)=[1-(\beta/t)^{k}]^{\gamma}$, $t>\theta$ and $\theta,k,\beta>0$. Thus the pdf and the cdf of $KwGMO-EEP$ distribution are given by

$$f^{KwGMOEEP}(t;a,b,\alpha,\theta,\beta,k,\gamma)=\frac{ab\theta\alpha^{\theta}\gamma k\theta^{k}t^{-(k+1)}[1-(\theta/t)^{k}]^{\gamma-1}[1-[1-(\theta/t)^{k}]^{\gamma}]^{\theta-1}}{[1-\bar{\alpha}[1-[1-(\theta/t)^{k}]^{\gamma}]]^{\theta+1}}$$

$$\left[1-\left\{\frac{\alpha[1-[1-(\theta/t)^{k}]^{\gamma}]}{1-\bar{\alpha}[1-[1-(\theta/t)^{k}]^{\gamma}]}\right\}^{\theta}\right]^{a-1}\left[1-\left[1-\left\{\frac{\alpha[1-[1-(\theta/t)^{k}]^{\gamma}]}{1-\bar{\alpha}[1-[1-(\theta/t)^{k}]^{\gamma}]}\right\}^{\theta}\right]^{a}\right]^{b-1}$$

cdf: $F^{KwGMOEEP}(t;a,b,\alpha,\theta,\beta,k,\gamma)=1-\left[1-\left[1-\left\{\frac{\alpha[1-[1-(\theta/t)^{k}]^{\gamma}]}{1-\bar{\alpha}[1-[1-(\theta/t)^{k}]^{\gamma}]}\right\}^{\theta}\right]^{a}\right]^{b}$

sf: $\bar{F}^{KwGMOEEP}(t;a,b,\alpha,\theta,\beta,k,\gamma)=\left[1-\left[1-\left\{\frac{\alpha[1-[1-(\theta/t)^{k}]^{\gamma}]}{1-\bar{\alpha}[1-[1-(\theta/t)^{k}]^{\gamma}]}\right\}^{\theta}\right]^{a}\right]^{b}$

hrf: $h^{KwGMOEEP}(t;a,b,\alpha,\theta,\beta,k,\gamma)=$

$$\frac{ab\theta\alpha^{\theta}\gamma k\theta^{k}t^{-(k+1)}[1-(\theta/t)^{k}]^{\gamma-1}[1-[1-(\theta/t)^{k}]^{\gamma}]^{\theta-1}}{[1-\bar{\alpha}[1-[1-(\theta/t)^{k}]^{\gamma}]]^{\theta+1}}$$

$$\left[1-\left\{\frac{\alpha[1-[1-(\theta/t)^{k}]^{\gamma}]}{1-\bar{\alpha}[1-[1-(\theta/t)^{k}]^{\gamma}]}\right\}^{\theta}\right]^{a-1}\left[1-\left[1-\left\{\frac{\alpha[1-[1-(\theta/t)^{k}]^{\gamma}]}{1-\bar{\alpha}[1-[1-(\theta/t)^{k}]^{\gamma}]}\right\}^{\theta}\right]^{a}\right]^{-1}$$

rhrf: $r^{KwGMOEEP}(t;a,b,\alpha,\theta,\beta,k,\gamma)=$

$$\frac{ab\theta\alpha^{\theta}\gamma k\theta^{k}t^{-(k+1)}[1-(\theta/t)^{k}]^{\gamma-1}[1-[1-(\theta/t)^{k}]^{\gamma}]^{\theta-1}}{[1-\bar{\alpha}[1-[1-(\theta/t)^{k}]^{\gamma}]]^{\theta+1}}$$

$$\left[1-\left\{\frac{\alpha[1-[1-(\theta/t)^{k}]^{\gamma}]}{1-\bar{\alpha}[1-[1-(\theta/t)^{k}]^{\gamma}]}\right\}^{\theta}\right]^{a-1}\left[1-\left[1-\left\{\frac{\alpha[1-[1-(\theta/t)^{k}]^{\gamma}]}{1-\bar{\alpha}[1-[1-(\theta/t)^{k}]^{\gamma}]}\right\}^{\theta}\right]^{a}\right]^{b-1}$$

$$\left[1-\left[1-\left[1-\left\{\frac{\alpha[1-[1-(\theta/t)^{k}]^{\gamma}]}{1-\bar{\alpha}[1-[1-(\theta/t)^{k}]^{\gamma}]}\right\}^{\theta}\right]^{a}\right]^{b}\right]^{-1}$$

chrf: $H^{KwGMOEEP}(t;a,b,\alpha,\theta,\beta,k,\gamma)=-b\log\left[1-\left[1-\left\{\frac{\alpha[1-[1-(\theta/t)^{k}]^{\gamma}]}{1-\bar{\alpha}[1-[1-(\theta/t)^{k}]^{\gamma}]}\right\}^{\theta}\right]^{a}\right]$

## 5. General results for the Kumaraswamy Genralized Marshall-Olkin-G ($KwGMO-G$) family of distributions

In this section we derive some general results for the proposed $KwGMO-G$ family.

### 5.1 Expansions of pdf and sf

By using binomial expansion in (6), we obtain

$$f^{KwGMOG}(t;a,b,\alpha,\theta)=\frac{ab\theta\alpha^{\theta}g(t)\bar{G}(t)^{\theta-1}}{[1-\bar{\alpha}\bar{G}(t)]^{\theta+1}}\left[1-\left\{\frac{\alpha\bar{G}(t)}{1-\bar{\alpha}\bar{G}(t)}\right\}^{\theta}\right]^{a-1}\left[1-\left[1-\left\{\frac{\alpha\bar{G}(t)}{1-\bar{\alpha}\bar{G}(t)}\right\}^{\theta}\right]^{a}\right]^{b-1}$$

$$=\frac{ab\theta\alpha^{\theta}g(t)\bar{G}(t)^{\theta-1}}{[1-\bar{\alpha}\bar{G}(t)]^{\theta+1}}\left[1-\left\{\frac{\alpha\bar{G}(t)}{1-\bar{\alpha}\bar{G}(t)}\right\}^{\theta}\right]^{a-1}\sum_{j=0}^{b-1}\binom{b-1}{j}(-1)^{j}\left[1-\left\{\frac{\alpha\bar{G}(t)}{1-\bar{\alpha}\bar{G}(t)}\right\}^{\theta}\right]^{aj}$$

$$=ab\,f^{GMO}(t;\alpha,\theta)\sum_{j=0}^{b-1}\binom{b-1}{j}(-1)^{j}[F^{GMO}(t;\alpha,\theta)]^{a(j+1)-1}$$

$$=f^{GMO}(t;\alpha,\theta)\sum_{j=0}^{b-1}A_{j}\,[F^{GMO}(t;\alpha,\theta)]^{a(j+1)-1} \qquad (10)$$

$$=\sum_{j=0}^{b-1}\frac{A_{j}}{a(j+1)}\frac{d}{dt}[F^{GMO}(t;\alpha,\theta)]^{a(j+1)}$$

$$=\sum_{j=0}^{b-1}A'_{j}\,\frac{d}{dt}[F^{GMO}(t;\alpha,\theta)]^{a(j+1)} \qquad (11)$$

Where, $A'_{j}=\frac{(-1)^{j}}{(j+1)}b\binom{b-1}{j}$ and $A_{j}=A'_{j}\,a(j+1)$

Alternatively, we can expand the pdf as

$$= f^{GMO}(t;\alpha,\theta)\sum_{j=0}^{b-1} A_j \, [1-\overline{F}^{GMO}(t;\alpha,\theta)]^{a(j+1)-1}$$

$$= f^{GMO}(t;\alpha,\theta)\sum_{j=0}^{b-1} A_j \sum_{k=0}^{a(j+1)-1}\binom{a(j+1)-1}{k}(-1)^k[\overline{F}^{GMO}(t;\alpha,\theta)]^k$$

$$= f^{GMO}(t;\alpha,\theta)\sum_{k=0}^{a(j+1)-1} B_k\,[\overline{F}^{GMO}(t;\alpha,\theta)]^k \tag{12}$$

$$= \sum_{k=0}^{a(j+1)-1}\frac{B_k}{k+1}\frac{d}{dt}[\overline{F}^{GMO}(t;\alpha,\theta)]^{k+1}$$

$$= \sum_{k=0}^{a(j+1)-1} B'_k\,\frac{d}{dt}[\overline{F}^{GMO}(t;\alpha,\theta)]^{k+1} \tag{13}$$

$$= \sum_{k=0}^{a(j+1)-1} B'_k\,\frac{d}{dt}[\overline{F}^{GMO}(t;\alpha,\theta(k+1))]$$

$$= \sum_{k=0}^{a(j+1)-1} B'_k\, f^{GMO}(t;\alpha,\theta(k+1)) \tag{14}$$

Where $B'_k = \dfrac{1}{k+1}\displaystyle\sum_{j=o}^{b-1} A_j\,(-1)^{k+1}\binom{a(j+1)-1}{k}$ and $B_k = -B'_k(k+1)$

Similarly an expansion for the survival function of $KwGMO-G$ can be derives as

$$\overline{F}^{KwGMOG}(t;a,b,\alpha,\theta) = \left[1-\left[1-\left\{\frac{\alpha\,\overline{G}(t)}{1-\overline{\alpha}\,\overline{G}(t)}\right\}^{\theta}\right]^{a}\right]^{b}$$

$$= \sum_{l=0}^{b}\binom{b}{l}(-1)^l\left[1-\left\{\frac{\alpha\,\overline{G}(t)}{1-\overline{\alpha}\,\overline{G}(t)}\right\}^{\theta}\right]^{a\,l}$$

$$= \sum_{l=0}^{b}\binom{b}{l}(-1)^l\,[F^{GMO}(t;\alpha,\theta)]^{al}$$

$$= \sum_{l=0}^{b} C_l\,[F^{GMO}(t;\alpha,\theta)]^{al} \tag{15}$$

Where, $C_l = \binom{b}{l}(-1)^l$

Again $\overline{F}^{KwGMOG}(t;a,b,\alpha,\theta) = \displaystyle\sum_{l=0}^{b} C_l\,[1-\overline{F}^{GMO}(t;\alpha,\theta)]^{al}$

$$= \sum_{l=0}^{b} C_l \sum_{m=0}^{al} \binom{al}{m} (-1)^m [\bar{F}^{GMO}(t;\alpha,\theta)]^m$$

$$= \sum_{m=0}^{al} \delta_m [\bar{F}^{GMO}(t;\alpha,\theta)]^m$$

$$= \sum_{m=0}^{al} \delta_m [\bar{F}^{GMO}(t;\alpha,\theta\, m)] \tag{16}$$

Where, $\delta_m = \sum_{l=0}^{b} C_l \binom{al}{m} (-1)^m$

**5.2 Order statistics**

Suppose $T_1, T_2, \ldots T_n$ is a random sample from any $KwGMO-G$ distribution. Let $T_{i:n}$ denote the $i^{th}$ order statistics. The pdf of $T_{i:n}$ can be expressed as

$$f_{i:n}(t) = \frac{n!}{(i-1)!\,(n-i)!} f^{KwGMOG}(t)[1-\bar{F}^{KwGMOG}(t)]^{i-1} \bar{F}^{KwGMOG}(t)^{n-i}$$

$$= \frac{n!}{(i-1)!\,(n-i)!} f^{KwGMOG}(t)\bar{F}^{KwGMOG}(t)^{n-i} \sum_{k=0}^{i-1} \binom{i-1}{k} [-\bar{F}^{KwGMOG}(t)]^k$$

$$= \frac{n!}{(i-1)!\,(n-i)!} f^{KwGMOG}(t) \sum_{k=0}^{i-1} (-1)^k \binom{i-1}{k} \bar{F}^{KwGMOG}(t)^{n+k-i}$$

Using the general expansion of the $KwGMO-G$ distribution pdf and sf we get the pdf of the $i^{th}$ order statistics for of the $KwGMO-G$ as

$$f_{i:n}(t) = \frac{n!}{(i-1)!\,(n-i)!} \left\{ f^{GMO}(t;\alpha,\theta) \sum_{p=0}^{a(j+1)-1} B_p [\bar{F}^{GMO}(t;\alpha,\theta)]^p \right\} \times \sum_{k=0}^{i-1} (-1)^k \binom{i-1}{k} \left\{ \sum_{m=0}^{al} \delta_m [\bar{F}^{GMO}(t;\alpha,\theta)]^m \right\}^{n+k-i}$$

Where $B_p$ and $\delta_m$ defined in section 5.1

$$= \frac{n!}{(i-1)!\,(n-i)!} \sum_{k=0}^{i-1} (-1)^k \binom{i-1}{k} \left\{ f^{GMO}(t;\alpha,\theta) \sum_{p=0}^{a(j+1)-1} B_p [\bar{F}^{GMO}(t;\alpha,\theta)]^p \right\} \times \left\{ \sum_{m=0}^{al} \delta_m [\bar{F}^{GMO}(t;\alpha,\theta)]^{m(n+k+i)} \right\}$$

$$= \frac{n!}{(i-1)!\,(n-i)!} \sum_{k=0}^{i-1} (-1)^k \binom{i-1}{k} \left\{ f^{GMO}(t;\alpha,\theta) \sum_{m=0}^{al} \sum_{p=0}^{a(j+1)-1} \delta_m B_p [\bar{F}^{GMO}(t;\alpha,\theta)]^{m(n+k+i)+p} \right\}$$

$$= f^{GMO}(t;\alpha,\theta) \sum_{m=0}^{al} \sum_{p=0}^{a(j+1)-1} \eta_{m,p} [\bar{F}^{GMO}(t;\alpha,\theta)]^{m(n+k+i)+p} \tag{17}$$

$$= -\sum_{m=0}^{al}\sum_{p=0}^{a(j+1)-1} \frac{\eta_{m,p}}{m(n+k+i)+p+1}\frac{d}{dt}[\overline{F}^{GMO}(t;\alpha,\theta)]^{m(n+k+i)+p+1}$$

$$= \sum_{m=0}^{al}\sum_{p=0}^{a(j+1)-1} \eta'_{m,p}\frac{d}{dt}[\overline{F}^{GMO}(t;\alpha,\theta)]^{m(n+k+i)+p+1}$$

$$= \sum_{m=0}^{al}\sum_{p=0}^{a(j+1)-1} \eta'_{m,p}\frac{d}{dt}[\overline{F}^{GMO}(t;\alpha,\theta(m(n+k+i)+p+1))]$$

$$= \sum_{m=0}^{al}\sum_{p=0}^{a(j+1)-1} \eta'_{m,p} f^{GMO}(t;\alpha,\theta(m(n+k+i)+p+1)) \qquad (18)$$

Where, $\eta_{m,p} = n\delta_m B_p \binom{n-1}{i-1}\sum_{k=0}^{i-1}\binom{i-1}{k}(-1)^k$ and $\eta'_{m,p} = -\eta_{m,p}/[m(n+k+i)+p+1]$

### 5.3 Moments

The probability weighted moments (PWMs), first proposed by Greenwood et al. (1979), are expectations of certain functions of a random variable whose mean exists. The $(p,q,r)^{th}$ PWM of $T$ is defined by

$$\Gamma_{p,q,r} = \int_{-\infty}^{\infty} t^p [F(t)]^q [1-F(t)]^r f(t)\,dt$$

From equations (10) and (12) the $s^{th}$ moment of $T$ can be written either as

$$= f^{GMO}(t;\alpha,\theta)\sum_{j=0}^{\infty} A_j\,[F^{GMO}(t;\alpha,\theta)]^{a(j+1)-1} \;\; = f^{GMO}(t;\alpha,\theta)\sum_{k=0}^{\infty} B_k [\overline{F}^{GMO}(t;\alpha,\theta)]^k$$

$$E(T^s) = \int_{-\infty}^{+\infty} t^s\, f^{GMO}(t;\alpha,\theta)\sum_{j=0}^{b-1} A_j\,[F^{GMO}(t;\alpha,\theta)]^{a(j+1)-1}\,dt$$

$$= \sum_{j=0}^{b-1} A_j \int_{-\infty}^{+\infty} t^s [1-[\alpha\overline{G}(t)/1-\overline{\alpha}\,\overline{G}(t)]^\theta]^{a(j+1)-1}[\theta\alpha^\theta g(t)\overline{G}(t)^{\theta-1}/[1-\overline{\alpha}\,\overline{G}(t)]^{\theta+1}]dt$$

$$= \sum_{j=0}^{b-1} A_j\,\Gamma_{s,a(j+1)-1,0}$$

or $E(T^s) = \int_{-\infty}^{+\infty} t^s\; f^{GMO}(t;\alpha,\theta) \sum_{k=0}^{a(j+1)-1} B_k [\overline{F}^{GMO}(t;\alpha,\theta)]^k \; dt$

$$= \sum_{k=0}^{a(j+1)-1} B_k \int_{-\infty}^{+\infty} t^s [[\alpha\overline{G}(t)/1-\overline{\alpha}\,\overline{G}(t)]^\theta]^k [\theta\alpha^\theta g(t)\overline{G}(t)^{\theta-1}/[1-\overline{\alpha}\,\overline{G}(t)]^{\theta+1}]dt$$

$$= \sum_{k=0}^{a(j+1)-1} B_k\,\Gamma_{s,0,k}$$

Where

$$\Gamma_{p,q,r} = \int_{-\infty}^{\infty} t^{p}\, \{1-[\alpha\,\overline{G}(t)/1-\overline{\alpha}\,\overline{G}(t)]^{\theta}\}^{q}\{[\alpha\,\overline{G}(t)/1-\overline{\alpha}\,\overline{G}(t)]^{\theta}\}^{r}[\theta\alpha^{\theta}\, g(t)\overline{G}(t)^{\theta-1}/[1-\overline{\alpha}\,\overline{G}(t)]^{\theta+1}]d$$

is the PWM of $\text{GMO}(\alpha,\theta)$ distribution.

Therefore the moments of the $\text{KwGMO-G}\,(a,b,\alpha,\theta)$ can be expresses in terms of the PWMs of $\text{GMO}(\alpha,\theta)$ (Jayakumar and Mathew, 2008). The PWM method can generally be used for estimating parameters quantiles of generalized distributions. These moments have low variance and no severe biases, and they compare favourably with estimators obtained by maximum likelihood.

Proceeding as above we can derive $s^{th}$ moment of the $i^{th}$ order statistic $T_{i:n}$, in a random sample of size $n$ from $KwGMO-G$ on using equations (17) as

$$E(T^{s}{}_{i,n}) = \sum_{m=0}^{al} \sum_{p=0}^{a(j+1)-1} \eta_{m,p}\, \Gamma_{s,0,m(n+k+i)+p}$$

Where $A_j$, $B_k$ and $\eta_{m,p}$ defined in section 5.1 and 5.2 respectively.

### 5.4 Moment generating function

The moment generating function of $KwGMO-G$ family can be easily expressed in terms of those of the exponentiated $GMO$ (Jayakumar and Mathew, 2008) distribution using the results of section 5.1. For example using equation (13) it can be seen that

$$M_T(s) = E[e^{sT}] = \int_{-\infty}^{\infty} e^{st} f(t)\, dt = \int_{-\infty}^{\infty} e^{st} \sum_{k=0}^{a(j+1)-1} B'_k \frac{d}{dt}[\overline{F}^{\,GMO}(t;\alpha,\theta)]^{k+1}\, dt$$

$$= \sum_{k=0}^{a(j+1)-1} B'_k \int_{-\infty}^{\infty} e^{st} \frac{d}{dt}[\overline{F}^{\,GMO}(t;\alpha,\theta)]^{k+1}\, dt = \sum_{k=0}^{a(j+1)-1} B'_k\, M_X(s)$$

Where $M_X(s)$ is the mgf of a $GMO$ (Jayakumar and Mathew, 2008) distribution.

### 5.5 Rényi Entropy

The entropy of a random variable is a measure of uncertainty variation and has been used in various situations in science and engineering. The Rényi entropy is defined by

$$I_R(\delta) = (1-\delta)^{-1} \log\left( \int_{-\infty}^{\infty} f(t)^{\delta}\, dt \right)$$

where $\delta > 0$ and $\delta \neq 1$ For furthers details, see Song (2001). Using binomial expansion in (6) we can write

$f^{\,KwGMOG}(t;a,b,\alpha,\theta)^{\delta}$

$$=\left[\frac{ab\theta\alpha^{\theta}\, g(t)\overline{G}(t)^{\theta-1}}{[1-\overline{\alpha}\,\overline{G}(t)]^{\theta+1}}\left[1-\left\{\frac{\alpha\,\overline{G}(t)}{1-\overline{\alpha}\,\overline{G}(t)}\right\}^{\theta}\right]^{a-1}\left[1-\left[1-\left\{\frac{\alpha\,\overline{G}(t)}{1-\overline{\alpha}\,\overline{G}(t)}\right\}^{\theta}\right]^{a}\right]^{b-1}\right]^{\delta}$$

$$=\{ab\, f^{GMO}(t;\alpha,\theta)\,[F^{GMO}(t;\alpha,\theta)]^{a-1}[1-[F^{GMO}(t;\alpha,\theta)]^{a}]^{b-1}\}^{\delta}$$

$$=\{ab\}^{\delta} f^{GMO}(t;\alpha,\theta)^{\delta}\;[F^{GMO}(t;\alpha,\theta)]^{\delta(a-1)}[1-[F^{GMO}(t;\alpha,\theta)]^{a}]^{\delta(b-1)}$$

$$=\{ab\}^{\delta} f^{GMO}(t;\alpha,\theta)^{\delta}\;[F^{GMO}(t;\alpha,\theta)]^{\delta(a-1)}\sum_{i=0}^{\delta(b-1)}\binom{\delta(b-1)}{i}(-1)^{i}[F^{GMO}(t;\alpha,\theta)]^{ai}$$

$$=\{ab\}^{\delta} f^{GMO}(t;\alpha,\theta)^{\delta}\;\sum_{i=0}^{\delta(b-1)}\binom{\delta(b-1)}{i}(-1)^{i}[F^{GMO}(t;\alpha,\theta)]^{a(i+\delta)-\delta}$$

$$=\{ab\}^{\delta} f^{GMO}(t;\alpha,\theta)^{\delta}\;\sum_{i=0}^{\delta(b-1)}\binom{\delta(b-1)}{i}(-1)^{i}[F^{GMO}(t;\alpha,\theta)]^{a(i+\delta)-\delta}$$

$$=f^{GMO}(t;\alpha,\theta)^{\delta}\;\sum_{i=0}^{\delta(b-1)}Z_{i}[F^{GMO}(t;\alpha,\theta)]^{a(i+\delta)-\delta}$$

Thus the Rényi entropy of $T$ can be obtained as

$$I_{R}(\delta)=(1-\delta)^{-1}\log\left(\int_{-\infty}^{\infty} f^{GMO}(t;\alpha,\theta)^{\delta}\;\sum_{i=0}^{\delta(b-1)}Z_{i}[F^{GMO}(t;\alpha,\theta)]^{a(i+\delta)-\delta}dt\right)$$

$$=(1-\delta)^{-1}\log\left(\sum_{i=0}^{\delta(b-1)}Z_{i}\int_{-\infty}^{\infty} f^{GMO}(t;\alpha,\theta)^{\delta}\;[F^{GMO}(t;\alpha,\theta)]^{a(i+\delta)-\delta}dt\right)$$

Where, $Z_{i}=\{ab\}^{\delta}\binom{\delta(b-1)}{i}(-1)^{i}$

## 5.6 Quantile function and random sample generation

We shall now present a formula for generating $KwGMO-G$ random variable by using inversion method by inverting the cdf or the survival function.

$$\overline{F}^{KwGMOG}(t;a,b,\alpha,\theta)=\left[1-\left[1-\left\{\frac{\alpha\,\overline{G}(t)}{1-\overline{\alpha}\,\overline{G}(t)}\right\}^{\theta}\right]^{a}\right]^{b}=\overline{F}^{KwGMO}(t)^{1/b}=1-\left[1-\left\{\frac{\alpha\,\overline{G}(t)}{1-\overline{\alpha}\,\overline{G}(t)}\right\}^{\theta}\right]^{a}$$

$$\Rightarrow[1-\overline{F}^{KwGMO}(t)^{1/b}]^{1/a}=1-\left\{\frac{\alpha\,\overline{G}(t)}{1-\overline{\alpha}\,\overline{G}(t)}\right\}^{\theta}\Rightarrow[1-\{1-\overline{F}^{KwGMO}(t)^{1/b}\}^{1/a}]^{1/\theta}=\frac{\alpha\,\overline{G}(t)}{1-\overline{\alpha}\,\overline{G}(t)}$$

$$\Rightarrow[1-\overline{\alpha}\,\overline{G}(t)][1-\{1-\overline{F}^{KwGMO}(t)^{1/b}\}^{1/a}]^{1/\theta}=\alpha\,\overline{G}(t)$$

$$\Rightarrow[1-\{1-\overline{F}^{KwGMO}(t)^{1/b}\}^{1/a}]^{1/\theta}=[\alpha+\overline{\alpha}[1-\{1-\overline{F}^{KwGMO}(t)^{1/b}\}^{1/a}]^{1/\theta}]\overline{G}(t)$$

$$\Rightarrow 1-G(t)=\frac{[1-\{1-\overline{F}^{KwGMO}(t)^{1/b}\}^{1/a}]^{1/\theta}}{\alpha+\overline{\alpha}[1-\{1-\overline{F}^{KwGMO}(t)^{1/b}\}^{1/a}]^{1/\theta}}$$

$$\Rightarrow G(t) = 1 - \frac{[1-\{1-\{1-F^{KwGMO}(t)\}^{1/b}\}^{1/a}]^{1/\theta}}{\alpha + \bar{\alpha}[1-\{1-\{1-F^{KwGMO}(t)\}^{1/b}\}^{1/a}]^{1/\theta}}$$

$$\Rightarrow t = G^{-1}\left(1 - \frac{[1-\{1-\{1-F^{KwGMO}(t)\}^{1/b}\}^{1/a}]^{1/\theta}}{\alpha + \bar{\alpha}[1-\{1-\{1-F^{KwGMO}(t)\}^{1/b}\}^{1/a}]^{1/\theta}}\right)$$

to generate an random variable from $KwGMO-G$ first generate a $u \sim U(0,1)$ then use

$$\Rightarrow t = G^{-1}\left(1 - \frac{[1-\{1-\{1-u\}^{1/b}\}^{1/a}]^{1/\theta}}{\alpha + \bar{\alpha}[1-\{1-\{1-u\}^{1/b}\}^{1/a}]^{1/\theta}}\right) \qquad (19)$$

The $p^{\text{th}}$ Quantile $t_p$ for $KwGMO-G$ can be easily obtained from (19) as

$$\Rightarrow t_p = G^{-1}\left(1 - \frac{[1-\{1-\{1-p\}^{1/b}\}^{1/a}]^{1/\theta}}{\alpha + \bar{\alpha}[1-\{1-\{1-p\}^{1/b}\}^{1/a}]^{1/\theta}}\right)$$

For example, let the base line distribution be exponential with parameter $\lambda > 0,$ having pdf and cdf as $g(t:\lambda) = \lambda e^{-\lambda t}, t > 0$ and $G(t:\lambda) = 1 - e^{-\lambda t},$ respectively. Therefore the $p^{\text{th}}$ Quantile $t_p$ of $KwGMO-E$ is given by

$$t_p = -\frac{1}{\lambda}\log\left[1 - \left[1 - \frac{[1-\{1-\{1-p\}^{1/b}\}^{1/a}]^{1/\theta}}{\alpha + \bar{\alpha}[1-\{1-\{1-p\}^{1/b}\}^{1/a}]^{1/\theta}}\right]\right]$$

### 5.7 Asymptotes

Here we investigate the asymptotic shapes of the proposed family following the methods followed in Alizadeh *et al.*, (2015).

**Proposition 2.** The asymptotes of equations (6), (7) and (8) as $t \to 0$ are given by

$$f(t) \sim \frac{\theta a b g(t)}{\alpha}\left[1 - \left\{\frac{\alpha \bar{G}(t)}{1-\bar{\alpha}\bar{G}(t)}\right\}^{\theta}\right]^{a-1} \quad \text{as } G(t) \to 0$$

$$F(t) \sim 0 \quad \text{as } G(t) \to 0$$

$$h(t) \sim \frac{\theta a b g(t)}{\alpha}\left[1 - \left\{\frac{\alpha \bar{G}(t)}{1-\bar{\alpha}\bar{G}(t)}\right\}^{\theta}\right]^{a-1} \quad \text{as } G(t) \to 0$$

**Proposition 3**. The asymptotes of equations (6), (7) and (8) as $t \to \infty$ are given by

$$f(t) \sim ab\theta\alpha^{\theta} g(t)\bar{G}(t)^{\theta-1}[1-[1-\{\alpha\bar{G}(t)\}^{\theta}]^{a}]^{b-1} \quad \text{as } t \to \infty$$

$$1 - F(t) \sim [1-[1-\{\alpha\bar{G}(t)\}^{\theta}]^{a}]^{b} \quad \text{as } t \to \infty$$

$$h(t) \sim ab\theta\alpha^{\theta} g(t)\bar{G}(t)^{\theta-1}[1-[1-\{\alpha\bar{G}(t)\}^{\theta}]^{a}]^{-1} \quad \text{as } t \to \infty$$

## 6. Estimation

### 6.1 Maximum likelihood method

The model parameters of the $KwGMO-G$ distribution can be estimated by maximum likelihood. Let $t=(t_1,t_2,...t_n)^T$ be a random sample of size $n$ from $KwGMO-G$ with parameter vector $\boldsymbol{\theta}=(\theta,\alpha,a,b,\boldsymbol{\beta}^T)^T$, where $\boldsymbol{\beta}=(\beta_1,\beta_2,...\beta_q)^T$ corresponds to the parameter vector of the baseline distribution *G*. Then the log-likelihood function for $\boldsymbol{\theta}$ is given by

$$\ell=\ell(\boldsymbol{\theta})=n\log(ab)+n\log\theta+n\theta\log\alpha+\sum_{i=0}^{n}\log[g(t_i,\boldsymbol{\beta})]+(\theta-1)\sum_{i=0}^{n}\log[\overline{G}(t_i,\boldsymbol{\beta})]$$

$$-(\theta+1)\sum_{i=0}^{n}\log[1-\overline{\alpha}\,\overline{G}(t_i,\boldsymbol{\beta})]+(a-1)\sum_{i=0}^{n}\log\left[1-\left\{\frac{\alpha\overline{G}(t_i,\boldsymbol{\beta})}{1-\overline{\alpha}\,\overline{G}(t_i,\boldsymbol{\beta})}\right\}^{\theta}\right]$$

$$+(b-1)\sum_{i=0}^{n}\log\left[1-\left[1-\left\{\frac{\alpha\overline{G}(t_i,\boldsymbol{\beta})}{1-\overline{\alpha}\,\overline{G}(t_i,\boldsymbol{\beta})}\right\}^{\theta}\right]^{a}\right] \qquad (20)$$

This log-likelihood function can not be solved analytically because of its complex form but it can be maximized numerically by employing global optimization methods available with software's like R, SAS, Mathematica or by solving the nonlinear likelihood equations obtained by differentiating (20).

By taking the partial derivatives of the log-likelihood function with respect to $\theta,\alpha,a,b$ and $\boldsymbol{\beta}$ we obtain the components of the score vector

$$U_{\boldsymbol{\theta}}=(U_\theta,U_\alpha,U_a,U_b,U_{\beta^T})^T$$

$$U_\theta=\frac{\partial\ell}{\partial\theta}=\frac{n}{\theta}+n\log\alpha+\sum_{i=0}^{n}\log[\overline{G}(t_i,\boldsymbol{\beta})]-\sum_{i=0}^{n}\log[1-\overline{\alpha}\,\overline{G}(t_i,\boldsymbol{\beta})]$$

$$+(1-a)\sum_{i=0}^{n}\frac{1}{1-\left\{\frac{\alpha\overline{G}(t_i,\boldsymbol{\beta})}{1-\overline{\alpha}\,\overline{G}(t_i,\boldsymbol{\beta})}\right\}^{\theta}}\left\{\frac{\alpha\overline{G}(t_i,\boldsymbol{\beta})}{1-\overline{\alpha}\,\overline{G}(t_i,\boldsymbol{\beta})}\right\}^{\theta}\log\left[\frac{\alpha\overline{G}(t_i,\boldsymbol{\beta})}{1-\overline{\alpha}\,\overline{G}(t_i,\boldsymbol{\beta})}\right]$$

$$+a(b-1)\sum_{i=0}^{n}\left[\frac{\left[1-\left\{\frac{\alpha\overline{G}(t_i,\boldsymbol{\beta})}{1-\overline{\alpha}\,\overline{G}(t_i,\boldsymbol{\beta})}\right\}^{\theta}\right]^{a-1}\left\{\frac{\alpha\overline{G}(t_i,\boldsymbol{\beta})}{1-\overline{\alpha}\,\overline{G}(t_i,\boldsymbol{\beta})}\right\}^{\theta}\log\left[\frac{\alpha\overline{G}(t_i,\boldsymbol{\beta})}{1-\overline{\alpha}\,\overline{G}(t_i,\boldsymbol{\beta})}\right]}{1-\left[1-\left\{\frac{\alpha\overline{G}(t_i,\boldsymbol{\beta})}{1-\overline{\alpha}\,\overline{G}(t_i,\boldsymbol{\beta})}\right\}^{\theta}\right]^{a}}\right]$$

$$U_{\alpha} = \frac{\partial \ell}{\partial \alpha} = \frac{n\theta}{\alpha} - (\theta+1)\sum_{i=0}^{n} \frac{\overline{G}(t_i,\boldsymbol{\beta})}{1-\overline{\alpha}\,\overline{G}(t_i,\boldsymbol{\beta})}$$
$$+\theta(a-1)\sum_{i=0}^{n} \frac{\alpha^{\theta-1}\overline{G}(t_i,\boldsymbol{\beta})^{\theta} G(t)}{[\{1-\overline{\alpha}\,\overline{G}(t_i,\boldsymbol{\beta})\}^{\theta} - \{\alpha\overline{G}(t_i,\boldsymbol{\beta})\}^{\theta}]\{1-\overline{\alpha}\,\overline{G}(t_i,\boldsymbol{\beta})\}}$$
$$+a\theta(b-1)\sum_{i=0}^{n} \frac{[\{1-\overline{\alpha}\,\overline{G}(t_i,\boldsymbol{\beta})\}^{\theta} - \{\alpha\overline{G}(t_i,\boldsymbol{\beta})\}^{\theta}]^{a-1}\alpha^{\theta-1}\overline{G}(t_i,\boldsymbol{\beta})^{\theta} G(t)}{[\{1-\overline{\alpha}\,\overline{G}(t_i,\boldsymbol{\beta})\}^{a\theta} - [\{1-\overline{\alpha}\,\overline{G}(t_i,\boldsymbol{\beta})\}^{\theta} - \{\alpha\overline{G}(t_i,\boldsymbol{\beta})\}^{\theta}]^{a}]\{1-\overline{\alpha}\,\overline{G}(t_i,\boldsymbol{\beta})\}}$$

$$U_{a} = \frac{\partial \ell}{\partial a} = \frac{n}{a} + \sum_{i=0}^{n} \log\left[1-\left\{\frac{\alpha\overline{G}(t_i,\boldsymbol{\beta})}{1-\overline{\alpha}\,\overline{G}(t_i,\boldsymbol{\beta})}\right\}^{\theta}\right]$$

$$+(b-1)\sum_{i=0}^{n} \frac{1}{1-\left[1-\left\{\frac{\alpha\overline{G}(t_i,\boldsymbol{\beta})}{1-\overline{\alpha}\,\overline{G}(t_i,\boldsymbol{\beta})}\right\}^{\theta}\right]^{a}} \left[1-\left\{\frac{\alpha\overline{G}(t_i,\boldsymbol{\beta})}{1-\overline{\alpha}\,\overline{G}(t_i,\boldsymbol{\beta})}\right\}^{\theta}\right]^{a} \log\left[1-\left\{\frac{\alpha\overline{G}(t_i,\boldsymbol{\beta})}{1-\overline{\alpha}\,\overline{G}(t_i,\boldsymbol{\beta})}\right\}^{\theta}\right]$$

$$U_{b} = \frac{\partial \ell}{\partial b} = \frac{n}{b} + \sum_{i=0}^{n} \log\left[1-\left[1-\left\{\frac{\alpha\overline{G}(t_i,\boldsymbol{\beta})}{1-\overline{\alpha}\,\overline{G}(t_i,\boldsymbol{\beta})}\right\}^{\theta}\right]^{a}\right]$$

$$U_{\beta} = \frac{\partial \ell}{\partial \boldsymbol{\beta}} = \sum_{i=0}^{r} \frac{g^{(\boldsymbol{\beta})}(t_i,\boldsymbol{\beta})}{g(t_i,\boldsymbol{\beta})} + (1-\theta)\sum_{i=0}^{r} \frac{G^{(\boldsymbol{\beta})}(t_i,\boldsymbol{\beta})}{\overline{G}(t_i,\boldsymbol{\beta})} - (\theta+1)\sum_{i=0}^{n} \frac{\overline{\alpha}\,G^{(\boldsymbol{\beta})}(t_i,\boldsymbol{\beta})}{1-\overline{\alpha}\,\overline{G}(t_i,\boldsymbol{\beta})}$$
$$+\theta(a-1)\sum_{i=0}^{n} \frac{\alpha^{\theta}\overline{G}(t_i,\boldsymbol{\beta})^{\theta-1} G^{(\boldsymbol{\beta})}(t_i,\boldsymbol{\beta})}{[\{1-\overline{\alpha}\,\overline{G}(t_i,\boldsymbol{\beta})\}^{\theta} - \{\alpha\overline{G}(t_i,\boldsymbol{\beta})\}^{\theta}]\{1-\overline{\alpha}\,\overline{G}(t_i,\boldsymbol{\beta})\}}$$
$$+a\theta(b-1)\sum_{i=0}^{n} \frac{[\{1-\overline{\alpha}\,\overline{G}(t_i,\boldsymbol{\beta})\}^{\theta} - \{\alpha\overline{G}(t_i,\boldsymbol{\beta})\}^{\theta}]^{a-1}\alpha^{\theta}\overline{G}(t_i,\boldsymbol{\beta})^{\theta-1} G^{(\boldsymbol{\beta})}(t_i,\boldsymbol{\beta})}{[\{1-\overline{\alpha}\,\overline{G}(t_i,\boldsymbol{\beta})\}^{a\theta} - [\{1-\overline{\alpha}\,\overline{G}(t_i,\boldsymbol{\beta})\}^{\theta} - \{\alpha\overline{G}(t_i,\boldsymbol{\beta})\}^{\theta}]^{a}]\{1-\overline{\alpha}\,\overline{G}(t_i,\boldsymbol{\beta})\}}$$

**6.2 Asymptotic standard error and confidence interval for the mles**

The asymptotic variance-covariance matrix of the MLEs of parameters can obtained by inverting the Fisher information matrix $\mathrm{I}(\boldsymbol{\theta})$ which can be derived using the second partial derivatives of the log-likelihood function with respect to each parameter. The $i\,j^{th}$ elements of $\mathrm{I}_n(\boldsymbol{\theta})$ are given by $\mathrm{I}_{ij} = -E\left(\frac{\partial^2 l(\boldsymbol{\theta})}{\partial\theta_i\,\partial\theta_j}\right), \quad i,\ j = 1,\ 2,\ \cdots, 3+q$

The exact evaluation of the above expectations may be cumbersome. In practice one can estimate $\mathrm{I}_n(\boldsymbol{\theta})$ by the observed Fisher's information matrix $\hat{\mathrm{I}}_n(\hat{\boldsymbol{\theta}})$ is defined as:

$$\hat{\mathrm{I}}_{ij} \approx \left(-\frac{\partial^2 l(\boldsymbol{\theta})}{\partial\theta_i\,\partial\theta_j}\right)_{\boldsymbol{\theta}=\hat{\boldsymbol{\theta}}}, \quad i,\ j = 1,\ 2,\cdots, 3+q$$

Using the general theory of MLEs under some regularity conditions on the parameters as $n \to \infty$ the asymptotic distribution of $\sqrt{n}\,(\hat{\boldsymbol{\theta}} - \boldsymbol{\theta})$ is $N_k(0, V_n)$ where $V_n = (v_{jj}) = \mathrm{I}_n^{-1}(\boldsymbol{\theta})$. The asymptotic behaviour remains valid if $V_n$ is replaced by $\hat{V}_n = \hat{\mathrm{I}}^{-1}(\hat{\boldsymbol{\theta}})$. This This result can be used to provide large sample standard errors and also construct confidence intervals for the model parameters. Thus an approximate standard error and $(1-\gamma/2)100\%$ confidence interval for the mle of $j^{\text{th}}$ parameter $\theta_j$ are respectively given by $\sqrt{\hat{v}_{jj}}$ and $\hat{\theta}_j \pm Z_{\gamma/2}\sqrt{\hat{v}_{jj}}$ , where $Z_{\gamma/2}$ is the $\gamma/2$ point of standard normal distribution.

As an illustration on the MLE method its large sample standard errors, confidence interval in the case of $KwGMO-E(a,b,\alpha,\theta,\lambda)$ is discussed in an appendix.

**6.3 Real life applications**

In this subsection, we consider fitting of four real data sets to compare the proposed $KwGMO-G$ distribution with another recently introduced model $GMOKw-G$ (Handique and Chakraborty, 2015) by taking as Weibull as the base line $G$ distribution. The mles of the parameters are by obtained by numerical maximization of log likelihood function. Standard errors (se) and 95% confidence intervals are also derived using large sample approach (see appendix). In order to compare the distributions, we have considered known criteria like AIC (Akaike Information Criterion), BIC (Bayesian Information Criterion), CAIC (Consistent Akaike Information Criterion) and HQIC (Hannan-Quinn Information Criterion). It may be noted that $AIC = 2k - 2l$; $BIC = k\log(n) - 2l$; $CAIC = AIC + 2k(k+1)/(n-k-1)$ and $HQIC = 2k\log[\log(n)] - 2l$. Where $k$ the number of parameters is, $n$ the sample size and $l$ is the maximized value of the log-likelihood function under the considered model. In these applications method of maximum likelihood will be used to obtain the estimate of parameters.

**Example I:**

Here we consider the following data set of 346 nicotine measurements made from several brands of cigarettes in 1998. The data have been collected by the Federal Trade Commission which is an independent agency of the US government, whose main mission is the promotion of consumer protection. [http: //www.ftc.gov/ reports/tobacco or http: // pw1.netcom.com/ rdavis2/ smoke. html.]

{1.3, 1.0, 1.2, 0.9, 1.1, 0.8, 0.5, 1.0, 0.7, 0.5, 1.7, 1.1, 0.8, 0.5, 1.2, 0.8, 1.1, 0.9, 1.2, 0.9, 0.8, 0.6, 0.3, 0.8, 0.6, 0.4, 1.1, 1.1, 0.2, 0.8, 0.5, 1.1, 0.1, 0.8, 1.7, 1.0, 0.8, 1.0, 0.8, 1.0, 0.2, 0.8, 0.4, 1.0, 0.2, 0.8, 1.4, 0.8, 0.5, 1.1, 0.9, 1.3, 0.9, 0.4, 1.4, 0.9, 0.5, 1.7, 0.9, 0.8, 0.8, 1.2, 0.9, 0.8, 0.5, 1.0, 0.6, 0.1, 0.2, 0.5, 0.1, 0.1, 0.9, 0.6, 0.9, 0.6, 1.2, 1.5, 1.1, 1.4, 1.2, 1.7, 1.4, 1.0,

0.7, 0.4, 0.9, 0.7, 0.8, 0.7, 0.4, 0.9, 0.6, 0.4, 1.2, 2.0, 0.7, 0.5, 0.9, 0.5, 0.9, 0.7, 0.9, 0.7, 0.4, 1.0, 0.7, 0.9, 0.7, 0.5, 1.3, 0.9, 0.8, 1.0, 0.7, 0.7, 0.6, 0.8, 1.1, 0.9, 0.9, 0.8, 0.8, 0.7, 0.7, 0.4, 0.5, 0.4, 0.9, 0.9, 0.7, 1.0, 1.0, 0.7, 1.3, 1.0, 1.1, 1.1, 0.9, 1.1, 0.8, 1.0, 0.7, 1.6, 0.8, 0.6, 0.8, 0.6, 1.2, 0.9, 0.6, 0.8, 1.0, 0.5, 0.8, 1.0, 1.1, 0.8, 0.8, 0.5, 1.1, 0.8, 0.9, 1.1, 0.8, 1.2, 1.1, 1.2, 1.1, 1.2, 0.2, 0.5, 0.7, 0.2, 0.5, 0.6, 0.1, 0.4, 0.6, 0.2, 0.5, 1.1, 0.8, 0.6, 1.1, 0.9, 0.6, 0.3, 0.9, 0.8, 0.8, 0.6, 0.4, 1.2, 1.3, 1.0, 0.6, 1.2, 0.9, 1.2, 0.9, 0.5, 0.8, 1.0, 0.7, 0.9, 1.0, 0.1, 0.2, 0.1, 0.1, 1.1, 1.0, 1.1, 0.7, 1.1, 0.7, 1.8, 1.2, 0.9, 1.7, 1.2, 1.3, 1.2, 0.9, 0.7, 0.7, 1.2, 1.0, 0.9, 1.6, 0.8, 0.8, 1.1, 1.1, 0.8, 0.6, 1.0, 0.8, 1.1, 0.8, 0.5, 1.5, 1.1, 0.8, 0.6, 1.1, 0.8, 1.1, 0.8, 1.5, 1.1, 0.8, 0.4, 1.0, 0.8, 1.4, 0.9, 0.9, 1.0, 0.9, 1.3, 0.8, 1.0, 0.5, 1.0, 0.7, 0.5, 1.4, 1.2, 0.9, 1.1, 0.9, 1.1, 1.0, 0.9, 1.2, 0.9, 1.2, 0.9, 0.5, 0.9, 0.7, 0.3, 1.0, 0.6, 1.0, 0.9, 1.0, 1.1, 0.8, 0.5, 1.1, 0.8, 1.2, 0.8, 0.5, 1.5, 1.5, 1.0, 0.8, 1.0, 0.5, 1.7, 0.3, 0.6, 0.6, 0.4, 0.5, 0.5, 0.7, 0.4, 0.5, 0.8, 0.5, 1.3, 0.9, 1.3, 0.9, 0.5, 1.2, 0.9, 1.1, 0.9, 0.5, 0.7, 0.5, 1.1, 1.1, 0.5, 0.8, 0.6, 1.2, 0.8, 0.4, 1.3, 0.8, 0.5, 1.2, 0.7, 0.5, 0.9, 1.3, 0.8, 1.2, 0.9}

**Table 1:** MLEs, standard error's and 95% confidence intervals and the AIC, BIC, CAIC and HQIC values for the nicotine measurements data

| Parameters | $GMOKw-W$ | $KwGMO-W$ |
|---|---|---|
| $\hat{a}$ | 0.765<br>(0.026)<br>(0.71, 0.82) | 0.771<br>(0.151)<br>(0.48, 1.07) |
| $\hat{b}$ | 2.139<br>(0.774)<br>(0.62, 3.66) | 0.398<br>(0.175)<br>(0.05, 0.74) |
| $\hat{\lambda}$ | 4.271<br>(0.018)<br>(4.24, 4.31) | 2.672<br>(0.926)<br>(0.86, 4.49) |
| $\hat{\beta}$ | 2.919<br>(0.026)<br>(2.87, 2.97) | 2.720<br>(0.409)<br>(1.92, 3.52) |
| $\hat{\alpha}$ | 1.097<br>(0.309)<br>(0.49, 1.70) | 1.765<br>(0.901)<br>(-0.00096, 3.53) |
| $\hat{\theta}$ | 0.114<br>(0.042)<br>(0.03, 0.19) | 1.157<br>(0.521)<br>(0.14, 2.17) |
| log-likelihood ($l_{\max}$) | -111.75 | **-109.73** |
| AIC | 235.50 | **231.46** |
| BIC | 258.58 | **254.54** |
| CAIC | 235.75 | **231.71** |
| HQIC | 244.69 | **240.66** |

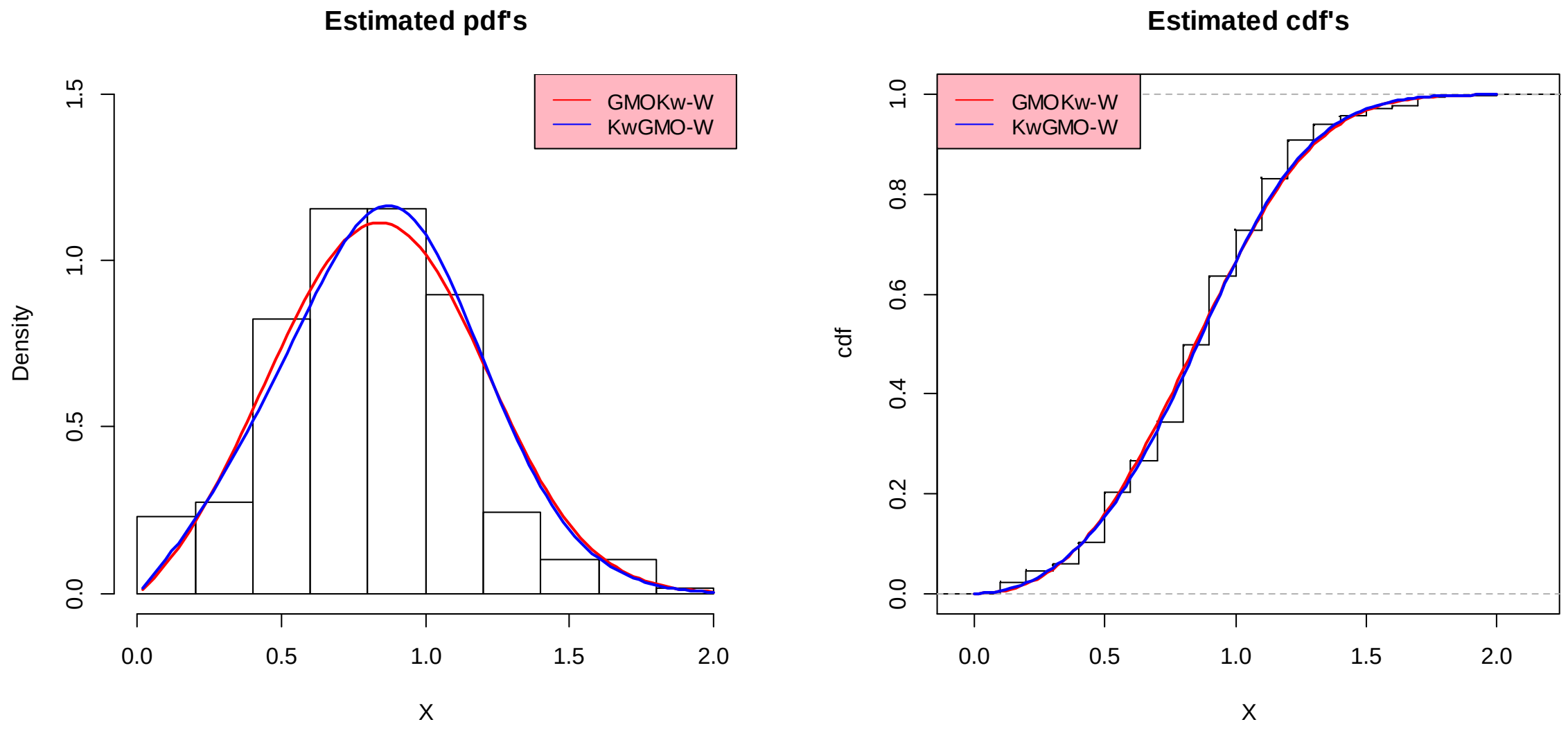

(a) (b)

**Fig: 3** Plots of the (a) observed histogram and estimated pdf's and (b) estimated cdf's for the $KwGMO-W$ and $GMOKw-W$ for example I.

**Example II:**

This data set consists of 100 observations of breaking stress of carbon fibres (in Gba) given by Nichols and Padgett (2006).

{0.98, 5.56, 5.08, 0.39, 1.57, 3.19, 4.90, 2.93, 2.85, 2.77, 2.76, 1.73, 2.48, 3.68, 1.08, 3.22, 3.75, 3.22, 3.70, 2.74, 2.73, 2.50, 3.60, 3.11, 3.27, 2.87, 1.47, 3.11, 4.42, 2.40, 3.15, 2.67,3.31, 2.81, 2.56, 2.17, 4.91, 1.59, 1.18, 2.48, 2.03, 1.69, 2.43, 3.39, 3.56, 2.83, 3.68, 2.00, 3.51, 0.85, 1.61, 3.28, 2.95, 2.81, 3.15, 1.92, 1.84, 1.22, 2.17, 1.61, 2.12, 3.09, 2.97, 4.20, 2.35, 1.41, 1.59, 1.12, 1.69, 2.79, 1.89, 1.87, 3.39, 3.33, 2.55, 3.68, 3.19, 1.71, 1.25, 4.70, 2.88, 2.96, 2.55, 2.59, 2.97, 1.57, 2.17, 4.38, 2.03, 2.82, 2.53, 3.31, 2.38, 1.36, 0.81, 1.17, 1.84, 1.80, 2.05, 3.65}.

**Table 2:** MLEs, standard errors and 95% confidence intervals (in parentheses) and the AIC, BIC, CAIC and HQIC values for the breaking stress of carbon fibres data.

| Parameters | $GMOKw-W$ | $KwGMO-W$ |
|---|---|---|
| $\hat{a}$ | 1.015<br>(0.071)<br>(0.88, 1.15) | 1.572<br>(1.584)<br>(-1.53, 4.68) |
| $\hat{b}$ | 0.385<br>(0.168)<br>(0.06, 0.71) | 1.029<br>(2.876)<br>(-4.61, 6.67) |
| $\hat{\lambda}$ | 0.803<br>(0.003)<br>(0.79, 0.81) | 0.575<br>(1.165)<br>(-1.71, 2.86) |
| $\hat{\beta}$ | 2.222<br>(0.004)<br>(2.21, 2.23) | 2.367<br>(1.128)<br>(0.16, 4.58) |
| $\hat{\alpha}$ | 1.482<br>(0.440)<br>(0.62, 2.34) | 0.620<br>(1.164)<br>(-1.66, 2.90) |
| $\hat{\theta}$ | 0.345<br>(0.164)<br>(0.02, 0.67) | 0.173<br>(0.375)<br>(-0.56, 0.91) |
| log-likelihood ($l_{\max}$) | -142.63 | **-141.29** |
| AIC | 297.26 | **294.58** |
| BIC | 312.89 | **310.21** |
| CAIC | 298.16 | **295.48** |
| HQIC | 303.59 | **300.92** |

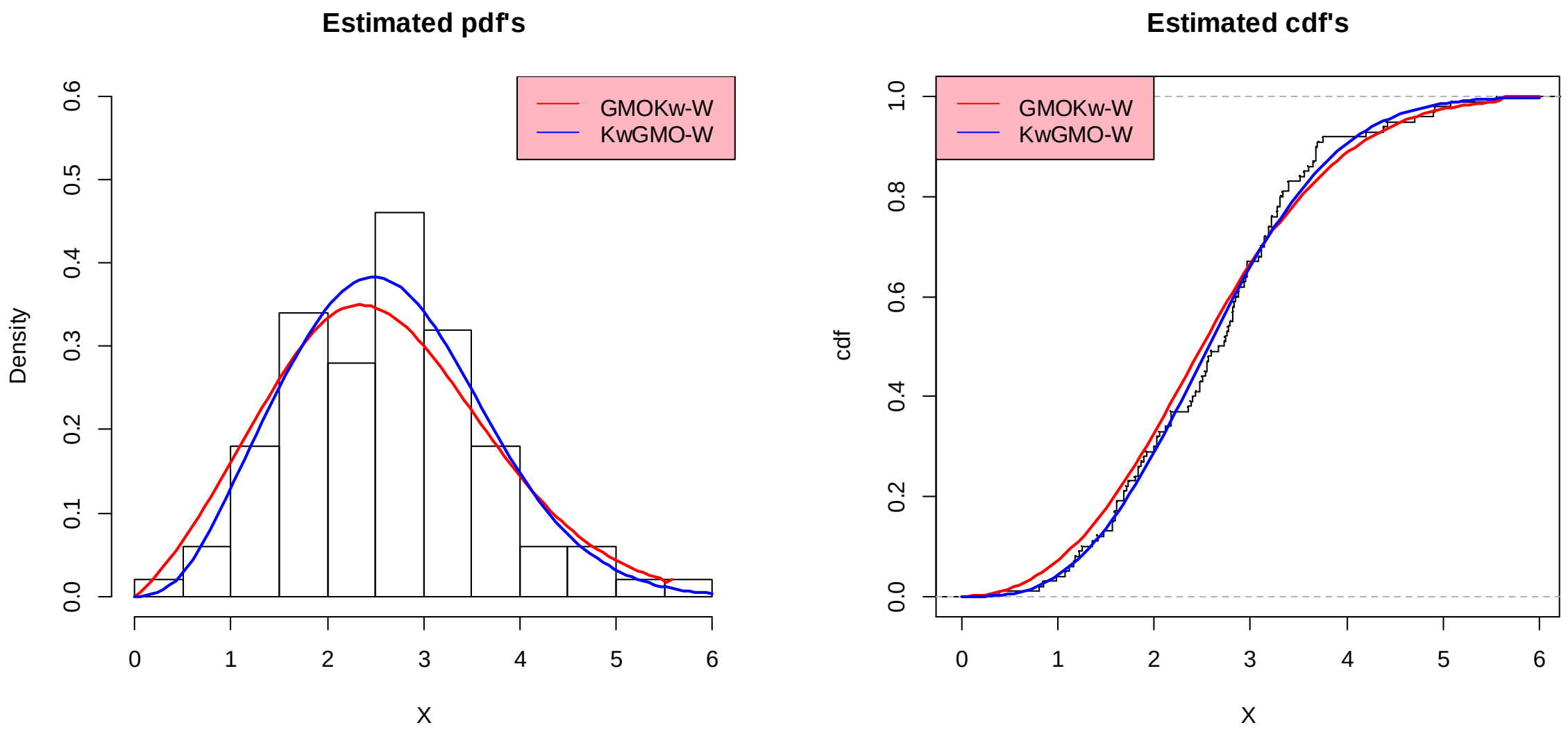

(a) (b)

**Fig: 4** Plots of the (a) observed histogram and estimated pdf's and (b) estimated cdf's for the $KwGMO-W$ and $GMOKw-W$ for example II.

**Example III:**

The Carbon Fibres data set of the 66 observations on the breaking stress of carbon fibres (in Gba) as reported in Nichols and Padgett (2006) is considered here.

{3.70, 2.74, 2.73, 2.50, 3.60, 3.11, 3.27, 2.87, 1.47, 3.11, 3.56, 4.42, 2.41, 3.19, 3.22, 1.69, 3.28, 3.09, 1.87, 3.15, 4.90, 1.57, 2.67, 2.93, 3.22, 3.39, 2.81, 4.20, 3.33, 2.55, 3.31, 3.31, 2.85, 1.25, 4.38, 1.84, 0.39, 3.68, 2.48, 0.85, 1.61, 2.79, 4.70, 2.03, 1.89, 2.88, 2.82, 2.05, 3.65, 3.75, 2.43, 2.95, 2.97, 3.39, 2.96, 2.35, 2.55, 2.59, 2.03, 1.61, 2.12, 3.15, 1.08, 2.56,1.80, 2.53.}

**Table 3:** MLEs, standard errors and 95% confidence intervals (in parentheses) and the AIC, BIC, CAIC and HQIC values for the breaking stress of carbon fibres data

| **Parameters** | $GMOKw-W$ | $KwGMO-W$ |
|---|---|---|
| $\hat{a}$ | 0.665<br>(0.339)<br>(0.00056, 1.33) | 0.592<br>(0.174)<br>(0.25, 0.93) |
| $\hat{b}$ | 7.809<br>(8.782)<br>(-9.40, 25.02) | 3.882<br>(2.307)<br>(-0.64, 8.40) |
| $\hat{\lambda}$ | 0.007<br>(0.004)<br>(-0.00084, 0.01) | 0.019<br>(0.015)<br>(-0.01, 0.05) |
| $\hat{\beta}$ | 3.369<br>(0.991)<br>(1.43, 5.31) | 4.489<br>(0.747)<br>(3.02, 5.95) |
| $\hat{\alpha}$ | 4.268<br>(4.572)<br>(-4.69, 13.23) | 2.552<br>(2.306)<br>(-1.97, 7.07) |
| $\hat{\theta}$ | 0.358<br>(0.232)<br>(-0.09, 0.81) | 0.037<br>(0.032)<br>(-0.025, 0.09) |
| log-likelihood ($l_{\max}$) | -85.11 | **-84.94** |
| AIC | 182.22 | **181.88** |
| BIC | 195.36 | **195.02** |
| CAIC | 183.64 | **183.30** |
| HQIC | 187.41 | **187.07** |

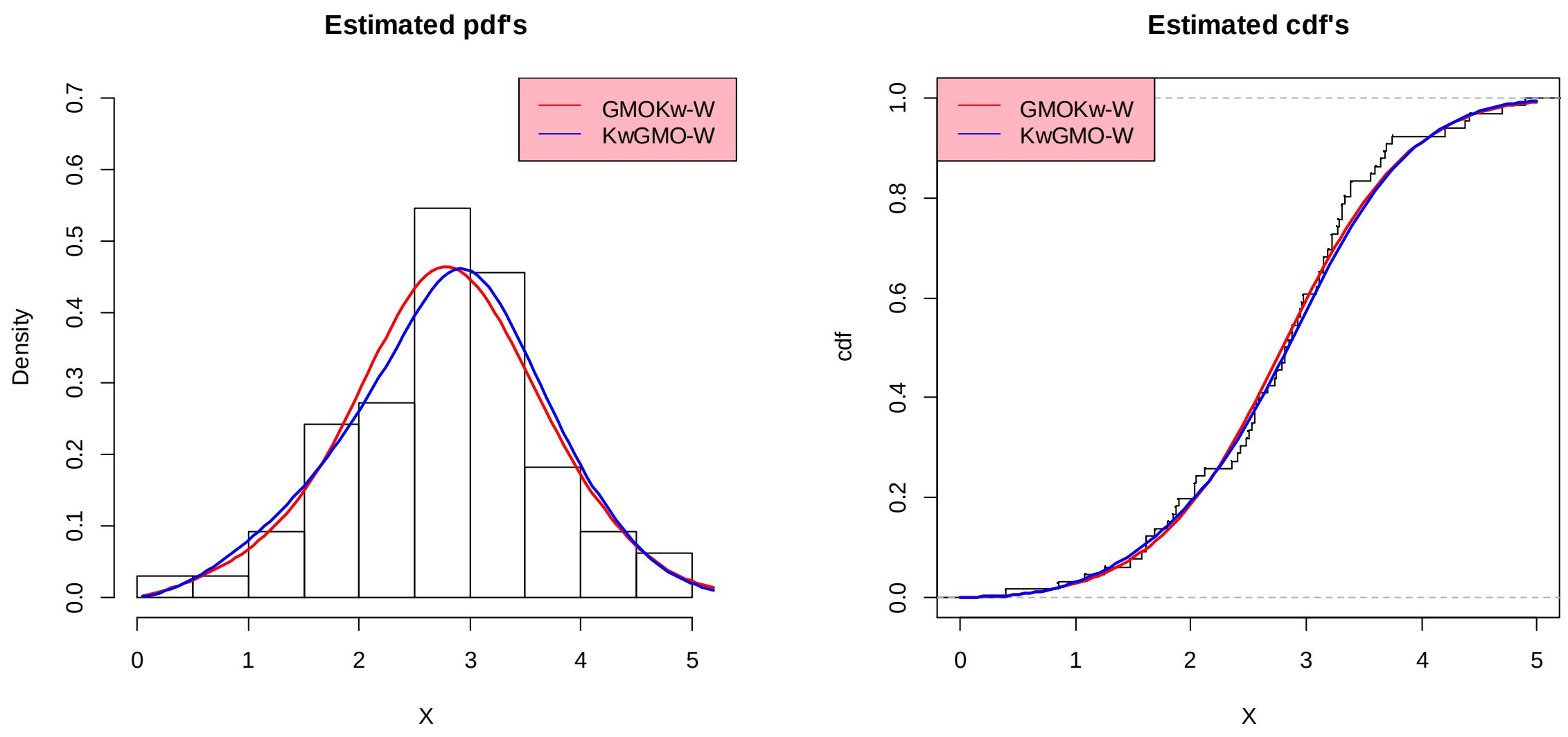

(a) (b)

**Fig: 5** Plots of the (a) observed histogram and estimated pdf's and (b) estimated cdf's for the $KwGMO-W$ and $GMOKw-W$ for example III.

**Example IV:**

Here, we use the following real data sets which gives the time to failure $(10^3 h)$ of turbocharger of one type of engine given in Xu et al. (2003).

{1.6, 2.0, 2.6, 3.0, 3.5, 3.9, 4.5, 4.6, 4.8, 5.0, 5.1, 5.3, 5.4, 5.6, 5.8, 6.0, 6.0, 6.1, 6.3, 6.5, 6.5, 6.7, 7.0, 7.1, 7.3, 7.3, 7.3, 7.7, 7.7, 7.8, 7.9, 8.0, 8.1, 8.3, 8.4, 8.4, 8.5, 8.7, 8.8, 9.0}

**Table 4:** MLEs, standard errors and 95% confidence intervals (in parentheses) and the AIC, BIC, CAIC and HQIC values for the data set

| Parameters | $GMOKw-W$ | $KwGMO-W$ |
|---|---|---|
| $\hat{a}$ | 1.178<br>**(0.017)**<br>(1.14, 1.21) | 1.156<br>**(0.602)**<br>(-0.02, 2.34) |
| $\hat{b}$ | 0.291<br>**(0.209)**<br>(-0.12, 0.70) | 1.785<br>**(2.746)**<br>(-3.59, 7.17) |
| $\hat{\lambda}$ | 0.617<br>**(0.002)**<br>(0.61, 0.62) | 0.007<br>**(0.005)**<br>(-0.0028, 0.02) |
| $\hat{\beta}$ | 1.855<br>**(0.003)**<br>(1.85, 1.86) | 4.206<br>**(0.727)**<br>(2.78, 5.63) |
| $\hat{\alpha}$ | 1.619<br>**(0.977)**<br>(-0.29, 3.53) | 0.047<br>**(0.087)**<br>(-0.12, 0.22) |
| $\hat{\theta}$ | 0.178<br>**(0.144)**<br>(-0.10, 0.46) | 0.002<br>**(0.023)**<br>(-0.04, 0.05) |
| log-likelihood($l_{\max}$) | -90.99 | **-80.33** |
| AIC | 193.98 | **172.66** |
| BIC | 204.11 | **182.79** |
| CAIC | 196.53 | **175.21** |
| HQIC | 197.65 | **176.33** |

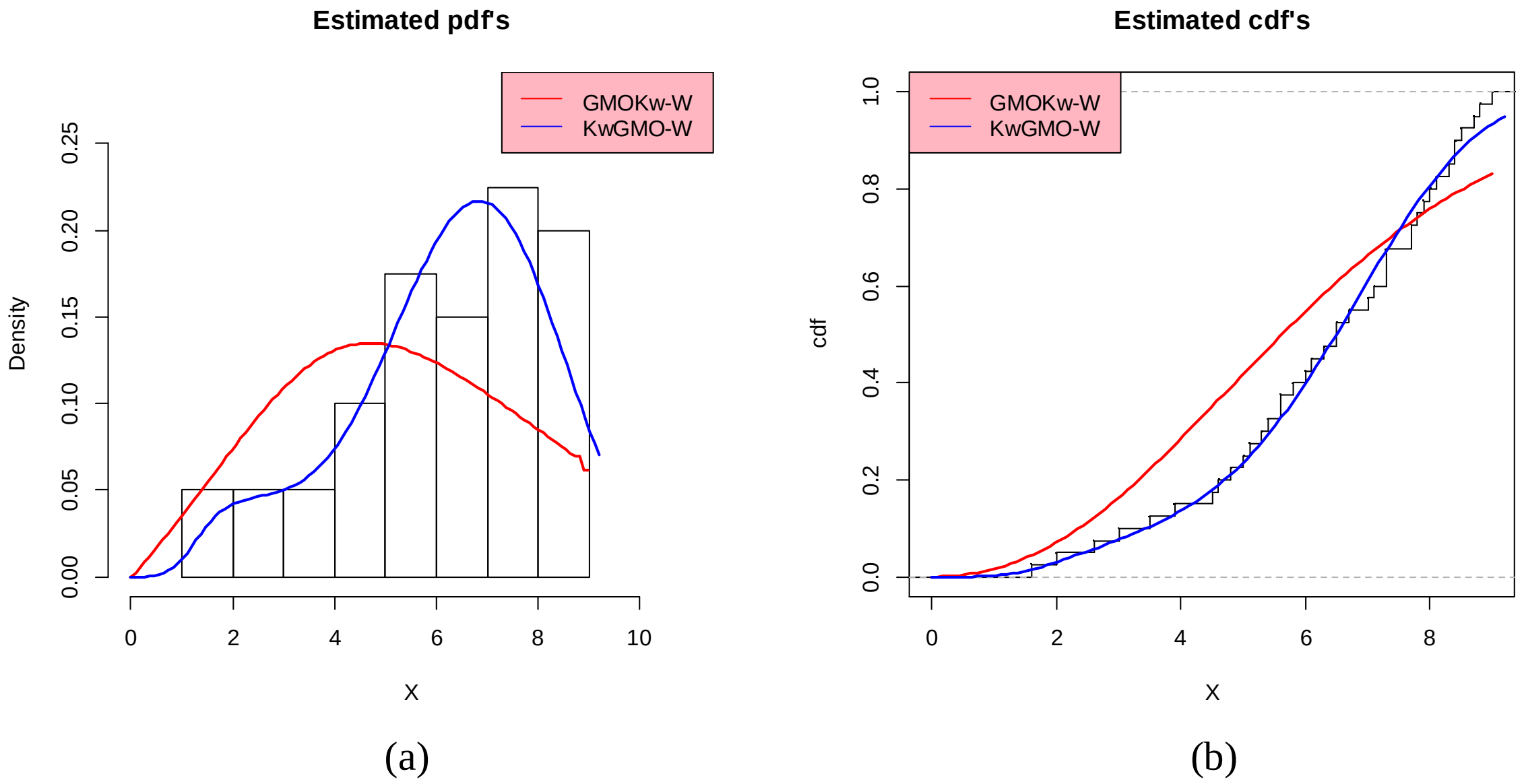

**Fig: 6** Plots of the (a) observed histogram and estimated pdf's and (b) estimated cdf's for the $KwGMO-W$ and $GMOKw-W$ for example IV.

In Tables 1, 2, 3 and 4 the MLEs, se's (in parentheses) and 95% confidence intervals (in parentheses) of the parameters for the fitted distributions along with AIC, BIC, CAIC and HQIC values are presented for example I, II, III and IV respectively. In all the four examples considered here based on the lowest values of the AIC, BIC, CAIC and HQIC, the $KwGMO-W$ distribution turns out to be a better distribution than $GMOKw-W$ distribution. A visual comparison of the closeness of the fitted densities with the observed histogram of the data and fitted cdfs with empirical cdfs for example I, II, III and IV are presented in the figures 3, 4, 5 and 6 respectively. These plots indicate that the proposed distributions provide a closer fit to these data.

## 7. Conclusion

Kumaraswamy Generalized Marshall-Olkin extended family of distributions is introduced and some of its important properties are studied. The maximum likelihood and moment method for estimating the parameters are also discussed. Four applications of real life data fitting shows good result in favour of the proposed family when compared to Generalized Marshall-Olkin Kumaraswamy extended family of distributions. As such the proposed family may be considered as a good contribution to the existing knowledge.

## Appendix: Maximum likelihood estimation for $KwGMO-E$

The pdf of the $KwGMO-E$ distribution is given by

$$f^{KwGMOE}(t;a,b,\alpha,\theta,\lambda)$$

$$=\frac{ab\theta\alpha^{\theta}\lambda e^{-\lambda t}[e^{-\lambda t}]^{\theta-1}}{[1-\overline{\alpha}e^{-\lambda t}]^{\theta+1}}\left[1-\left\{\frac{\alpha e^{-\lambda t}}{1-\overline{\alpha}e^{-\lambda t}}\right\}^{\theta}\right]^{a-1}\left[1-\left[1-\left\{\frac{\alpha e^{-\lambda t}}{1-\overline{\alpha}e^{-\lambda t}}\right\}^{\theta}\right]^{a}\right]^{b-1}$$

for $a>0, b>0, \alpha>0, \theta>0, \lambda>0, t>0$

For a random sample of size $n$ from this distribution, the log-likelihood function for the parameter vector $\boldsymbol{\theta}=(\theta,\alpha,a,b,\lambda)^T$ is given by

$$\ell=\ell(\boldsymbol{\theta})=n\log\theta+n\theta\log\alpha+n\log(ab)+n\log\lambda-(\theta-1)\lambda\sum_{i=0}^{n}t_i-(\theta+1)\sum_{i=0}^{n}\log[1-\overline{\alpha}e^{-\lambda t_i}]$$

$$+(a-1)\sum_{i=0}^{n}\log\left[1-\left\{\frac{\alpha e^{-\lambda t_i}}{1-\overline{\alpha}e^{-\lambda t_i}}\right\}^{\theta}\right]+(b-1)\sum_{i=0}^{n}\log\left[1-\left[1-\left\{\frac{\alpha e^{-\lambda t_i}}{1-\overline{\alpha}e^{-\lambda t_i}}\right\}^{\theta}\right]^{a}\right]$$

The components of the score vector $\boldsymbol{\theta}=(\theta,\alpha,a,b,\lambda)^T$ are

$$\frac{\partial\ell(\boldsymbol{\theta})}{\partial\theta}=\frac{n}{\theta}+n\log\alpha-\lambda\sum_{i=0}^{n}t_i-\sum_{i=0}^{n}\log[1-\overline{\alpha}e^{-\lambda t}]$$

$$+(1-a)\sum_{i=0}^{n}\frac{1}{1-\left\{\frac{\alpha e^{-\lambda t}}{1-\overline{\alpha}e^{-\lambda t}}\right\}^{\theta}}\left\{\frac{\alpha e^{-\lambda t}}{1-\overline{\alpha}e^{-\lambda t}}\right\}^{\theta}\log\left[\frac{\alpha e^{-\lambda t}}{1-\overline{\alpha}e^{-\lambda t}}\right]$$

$$+a(b-1)\sum_{i=0}^{n}\frac{1}{1-\left[1-\left\{\frac{\alpha e^{-\lambda t}}{1-\overline{\alpha}e^{-\lambda t}}\right\}^{\theta}\right]^{a}}\left[1-\left\{\frac{\alpha e^{-\lambda t}}{1-\overline{\alpha}e^{-\lambda t}}\right\}^{\theta}\right]^{a-1}\left\{\frac{\alpha e^{-\lambda t}}{1-\overline{\alpha}e^{-\lambda t}}\right\}^{\theta}\log\left[\frac{\alpha e^{-\lambda t}}{1-\overline{\alpha}e^{-\lambda t}}\right]$$

$$\frac{\partial\ell(\boldsymbol{\theta})}{\partial\alpha}=\frac{n\theta}{\alpha}-(\theta+1)\sum_{i=0}^{n}\frac{e^{-\lambda t_i}}{1-\overline{\alpha}e^{-\lambda t_i}}$$

$$+\theta(1-a)\sum_{i=0}^{n}\frac{(\alpha e^{-\lambda t_i}/1-\overline{\alpha}e^{-\lambda t_i})^{\theta-1}\left\{-\alpha e^{-2\lambda t_i}/(1-\overline{\alpha}e^{-\lambda t_i})^2+(e^{-\lambda t_i}/1-\overline{\alpha}e^{-\lambda t_i})\right\}}{1-(\alpha e^{-\lambda t_i}/1-\overline{\alpha}e^{-\lambda t_i})^{\theta}}$$

$$+a\theta(b-1)\sum_{i=0}^{n}\frac{(\alpha e^{-\lambda t_i}/1-\overline{\alpha}e^{-\lambda t_i})^{\theta-1}\{1-(\alpha e^{-\lambda t_i}/1-\overline{\alpha}e^{-\lambda t_i})^{\theta}\}^{a-1}}{1-\{1-(\alpha e^{-\lambda t_i}/1-\overline{\alpha}e^{-\lambda t_i})^{\theta}\}^{a}}$$

$$\left\{-\alpha e^{-2\lambda t_i}/(1-\overline{\alpha}e^{-\lambda t_i})^2+(e^{-\lambda t_i}/1-\overline{\alpha}e^{-\lambda t_i})\right\}$$

$$\frac{\partial\ell(\boldsymbol{\theta})}{\partial a}=\frac{n}{a}+\sum_{i=0}^{n}\log\left[1-\left\{\frac{\alpha e^{-\lambda t_i}}{1-\overline{\alpha}e^{-\lambda t_i}}\right\}^{\theta}\right]$$

$$+(1-b)\sum_{i=0}^{n}\frac{1}{1-\left[1-\left\{\dfrac{\alpha\,e^{-\lambda t_i}}{1-\overline{\alpha}\,e^{-\lambda t_i}}\right\}^{\theta}\right]^{a}}\left[1-\left\{\frac{\alpha\,e^{-\lambda t_i}}{1-\overline{\alpha}\,e^{-\lambda t_i}}\right\}^{\theta}\right]^{a}\log\left[1-\left\{\frac{\alpha\,e^{-\lambda t_i}}{1-\overline{\alpha}\,e^{-\lambda t_i}}\right\}^{\theta}\right]$$

$$\frac{\partial\ell(\boldsymbol{\theta})}{\partial b}=\frac{n}{b}+\sum_{i=0}^{n}\log\left[1-\left[1-\left\{\frac{\alpha e^{-\lambda t_i}}{1-\overline{\alpha}\,e^{-\lambda t_i}}\right\}^{\theta}\right]^{a}\right]$$

$$\frac{\partial\ell(\boldsymbol{\theta})}{\partial\lambda}=\frac{n}{\lambda}-(\theta-1)\sum_{i=0}^{n}t_i-(\theta+1)\sum_{i=0}^{n}\frac{\overline{\alpha}\,e^{-\lambda t_i}\,t_i}{1-\overline{\alpha}\,e^{-\lambda t_i}}$$

$$+\theta(1-a)\sum_{i=0}^{n}\frac{(\alpha e^{-\lambda t_i}/1-\overline{\alpha}\,e^{-\lambda t_i})^{\theta-1}\left\{-\alpha\overline{\alpha}\,e^{-2\lambda t_i}\,t_i/(1-\overline{\alpha}\,e^{-\lambda t_i})^{2}-(\alpha e^{-\lambda t_i}\,t_i/1-\overline{\alpha}\,e^{-\lambda t_i})\right\}}{1-(\alpha e^{-\lambda t_i}/1-\overline{\alpha}\,e^{-\lambda t_i})^{\theta}}$$

$$+a\theta(b-1)\sum_{i=0}^{n}\frac{(\alpha e^{-\lambda t_i}/1-\overline{\alpha}\,e^{-\lambda t_i})^{\theta-1}\{1-(\alpha e^{-\lambda t_i}/1-\overline{\alpha}\,e^{-\lambda t_i})^{\theta}\}^{a-1}}{1-\{1-(\alpha e^{-\lambda t_i}/1-\overline{\alpha}\,e^{-\lambda t_i})^{\theta}\}^{a}}$$

$$\times\left\{-\alpha\overline{\alpha}\,e^{-2\lambda t_i}\,t_i/(1-\overline{\alpha}\,e^{-\lambda t_i})^{2}-(\alpha e^{-\lambda t_i}\,t_i/1-\overline{\alpha}\,e^{-\lambda t_i})\right\}$$

The asymptotic variance covariance matrix for mles of the parameters of $KwGMO-E$ $(\theta,\alpha,a,b,\lambda)$ distribution is estimated by

$$\hat{\mathrm{I}}_{\mathrm{n}}^{\ -1}(\hat{\boldsymbol{\theta}})=\begin{pmatrix}\mathrm{var}(\hat{\theta}) & \mathrm{cov}(\hat{\theta},\hat{\alpha}) & \mathrm{cov}(\hat{\theta},\hat{a}) & \mathrm{cov}(\hat{\theta},\hat{b}) & \mathrm{cov}(\hat{\theta},\hat{\lambda})\\ \mathrm{cov}(\hat{\alpha},\hat{\theta}) & \mathrm{var}(\hat{\alpha}) & \mathrm{cov}(\hat{\alpha},\hat{a}) & \mathrm{cov}(\hat{\alpha},\hat{b}) & \mathrm{cov}(\hat{\alpha},\hat{\lambda})\\ \mathrm{cov}(\hat{a},\hat{\theta}) & \mathrm{cov}(\hat{a},\hat{\alpha}) & \mathrm{var}(\hat{a}) & \mathrm{cov}(\hat{a},\hat{b}) & \mathrm{cov}(\hat{a},\hat{\lambda})\\ \mathrm{cov}(\hat{b},\hat{\theta}) & \mathrm{cov}(\hat{b},\hat{\alpha}) & \mathrm{cov}(\hat{b},\hat{a}) & \mathrm{var}(\hat{b}) & \mathrm{cov}(\hat{b},\hat{\lambda})\\ \mathrm{cov}(\hat{\lambda},\hat{\theta}) & \mathrm{cov}(\hat{\lambda},\hat{\alpha}) & \mathrm{cov}(\hat{\lambda},\hat{a}) & \mathrm{cov}(\hat{\lambda},\hat{b}) & \mathrm{var}(\hat{\lambda})\end{pmatrix}$$

Where the elements of the information matrix $\hat{\mathbf{I}}_{\mathrm{n}}\,(\hat{\boldsymbol{\theta}})=\left(-\left(\dfrac{\partial^{2}l(\boldsymbol{\theta})}{\partial\theta_i\,\partial\theta_j}\right)\right)_{\boldsymbol{\theta}=\hat{\boldsymbol{\theta}}}$ can be derived using the following second partial derivatives:

$$\frac{\partial^{2}\ell}{\partial\theta^{2}}=\frac{n}{\theta^{2}}+(1-a)\sum_{i=0}^{n}\frac{(\alpha e^{-\lambda t}/1-\overline{\alpha}\,e^{-\lambda t})^{2\theta}\log(\alpha e^{-\lambda t}/1-\overline{\alpha}\,e^{-\lambda t})^{2}}{\{1-(\alpha e^{-\lambda t}/1-\overline{\alpha}\,e^{-\lambda t})^{\theta}\}^{2}}$$

$$+(1-a)\sum_{i=0}^{n}\frac{(\alpha e^{-\lambda t}/1-\overline{\alpha}\,e^{-\lambda t})^{\theta}\log(\alpha e^{-\lambda t}/1-\overline{\alpha}\,e^{-\lambda t})^{2}}{1-(\alpha e^{-\lambda t}/1-\overline{\alpha}\,e^{-\lambda t})^{\theta}}$$

$$+(1-b)\sum_{i=0}^{n}\frac{a^2\{1-(\alpha e^{-\lambda t}/1-\overline{\alpha}\, e^{-\lambda t})^{\theta}\}^{2(a-1)}(\alpha e^{-\lambda t}/1-\overline{\alpha}\, e^{-\lambda t})^{2\theta}\log(\alpha e^{-\lambda t}/1-\overline{\alpha}\, e^{-\lambda t})^{2}}{[1-\{1-(\alpha e^{-\lambda t}/1-\overline{\alpha}\, e^{-\lambda t})^{\theta}\}]^{2}}$$

$$+(1-b)\sum_{i=0}^{n}\frac{a(a-1)\{1-(\alpha e^{-\lambda t}/1-\overline{\alpha}\, e^{-\lambda t})^{\theta}\}^{a-2}(\alpha e^{-\lambda t}/1-\overline{\alpha}\, e^{-\lambda t})^{2\theta}\log(\alpha e^{-\lambda t}/1-\overline{\alpha}\, e^{-\lambda t})^{2}}{1-\{1-(\alpha e^{-\lambda t}/1-\overline{\alpha}\, e^{-\lambda t})^{\theta}\}}$$

$$-(1-b)\sum_{i=0}^{n}\frac{a\,\{1-(\alpha e^{-\lambda t}/1-\overline{\alpha}\, e^{-\lambda t})^{\theta}\}^{a-1}(\alpha e^{-\lambda t}/1-\overline{\alpha}\, e^{-\lambda t})^{\theta}\log(\alpha e^{-\lambda t}/1-\overline{\alpha}\, e^{-\lambda t})^{2}}{1-\{1-(\alpha e^{-\lambda t}/1-\overline{\alpha}\, e^{-\lambda t})^{\theta}\}}$$

$$\frac{\partial^2\ell}{\partial\alpha^2}=-\frac{n\theta}{\alpha^2}+(\theta+1)\sum_{i=0}^{n}\frac{e^{-2\lambda t_i}}{(1-\overline{\alpha}\, e^{-\lambda t_i})^2}$$

$$+\theta(1-a)\sum_{i=0}^{n}\frac{(\alpha e^{-\lambda t_i}/1-\overline{\alpha}\, e^{-\lambda t_i})^{\theta-1}\left\{2\alpha e^{-3\lambda t_i}/(1-\overline{\alpha}\, e^{-\lambda t_i})^{3}-(2e^{-2\lambda t_i}/(1-\overline{\alpha}\, e^{-\lambda t_i})^{2}\right\}}{1-(\alpha e^{-\lambda t_i}/1-\overline{\alpha}\, e^{-\lambda t_i})^{\theta}}$$

$$+\theta^2(1-a)\sum_{i=0}^{n}\frac{(\alpha e^{-\lambda t_i}/1-\overline{\alpha}\, e^{-\lambda t_i})^{2(\theta-1)}\left\{-\alpha e^{-2\lambda t_i}/(1-\overline{\alpha}\, e^{-\lambda t_i})^{2}+(e^{-\lambda t_i}/(1-\overline{\alpha}\, e^{-\lambda t_i}))\right\}^{2}}{\{1-(\alpha e^{-\lambda t_i}/1-\overline{\alpha}\, e^{-\lambda t_i})^{\theta}\}^{2}}$$

$$+\theta(\theta-1)(1-a)\sum_{i=0}^{n}\frac{(\alpha e^{-\lambda t_i}/1-\overline{\alpha}\, e^{-\lambda t_i})^{\theta-2}\left\{-\alpha e^{-2\lambda t_i}/(1-\overline{\alpha}\, e^{-\lambda t_i})^{2}+(e^{-\lambda t_i}/(1-\overline{\alpha}\, e^{-\lambda t_i}))\right\}^{2}}{1-(\alpha e^{-\lambda t_i}/1-\overline{\alpha}\, e^{-\lambda t_i})^{\theta}}$$

$$+a^2\theta^2(1-b)\sum_{i=0}^{n}\frac{(\alpha e^{-\lambda t_i}/1-\overline{\alpha}\, e^{-\lambda t_i})^{2(\theta-1)}\{1-(\alpha e^{-\lambda t_i}/1-\overline{\alpha}\, e^{-\lambda t_i})^{\theta}\}^{2(a-1)}}{\{1-\{1-(\alpha e^{-\lambda t_i}/1-\overline{\alpha}\, e^{-\lambda t_i})^{\theta}\}^{a}\}^{2}}$$

$$\times\left\{-\alpha e^{-2\lambda t_i}/(1-\overline{\alpha}\, e^{-\lambda t_i})^{2}+(e^{-\lambda t_i}/1-\overline{\alpha}\, e^{-\lambda t_i})\right\}^{2}$$

$$-a\theta(1-b)\sum_{i=0}^{n}\frac{(\alpha e^{-\lambda t_i}/1-\overline{\alpha}\, e^{-\lambda t_i})^{\theta-1}\{1-(\alpha e^{-\lambda t_i}/1-\overline{\alpha}\, e^{-\lambda t_i})^{\theta}\}^{a-1}}{1-\{1-(\alpha e^{-\lambda t_i}/1-\overline{\alpha}\, e^{-\lambda t_i})^{\theta}\}^{a}}$$

$$\times\left\{2\alpha e^{-3\lambda t_i}/(1-\overline{\alpha}\, e^{-\lambda t_i})^{2}-(2e^{-2\lambda t_i}/(1-\overline{\alpha}\, e^{-\lambda t_i})^{2})\right\}$$

$$+a\,(a-1)\theta^2(1-b)\sum_{i=0}^{n}\frac{(\alpha e^{-\lambda t_i}/1-\overline{\alpha}\, e^{-\lambda t_i})^{2(\theta-1)}\{1-(\alpha e^{-\lambda t_i}/1-\overline{\alpha}\, e^{-\lambda t_i})^{\theta}\}^{a-2}}{1-\{1-(\alpha e^{-\lambda t_i}/1-\overline{\alpha}\, e^{-\lambda t_i})^{\theta}\}^{a}}$$

$$\times\left\{-\alpha e^{-2\lambda t_i}/(1-\overline{\alpha}\, e^{-\lambda t_i})^{2}+(e^{-\lambda t_i}/1-\overline{\alpha}\, e^{-\lambda t_i})\right\}^{2}$$

$$+a\;\theta(\theta-1)(1-b)\sum_{i=0}^{n}\frac{(\alpha e^{-\lambda t_i}/1-\overline{\alpha}\, e^{-\lambda t_i})^{\theta-2}\{1-(\alpha e^{-\lambda t_i}/1-\overline{\alpha}\, e^{-\lambda t_i})^{\theta}\}^{a-1}}{1-\{1-(\alpha e^{-\lambda t_i}/1-\overline{\alpha}\, e^{-\lambda t_i})^{\theta}\}^{a}}$$

$$\times\left\{-\alpha e^{-2\lambda t_i}/(1-\overline{\alpha}\, e^{-\lambda t_i})^{2}+(e^{-\lambda t_i}/1-\overline{\alpha}\, e^{-\lambda t_i})\right\}^{2}$$

$$\frac{\partial^2\ell}{\partial a^2}=-\frac{n}{a^2}+(1-b)\sum_{i=0}^{n}\frac{\{1-(\alpha e^{-\lambda t_i}/1-\overline{\alpha}\, e^{-\lambda t_i})^{\theta}\}^{2a}\log\{1-(\alpha e^{-\lambda t_i}/1-\overline{\alpha}\, e^{-\lambda t_i})^{\theta}\}^{2}}{\{1-\{1-(\alpha e^{-\lambda t_i}/1-\overline{\alpha}\, e^{-\lambda t_i})^{\theta}\}^{a}\}^{2}}$$

$$+(1-b)\sum_{i=0}^{n}\frac{\{1-(\alpha e^{-\lambda t_i}/1-\overline{\alpha}\, e^{-\lambda t_i})^{\theta}\}^{a}\log\{1-(\alpha e^{-\lambda t_i}/1-\overline{\alpha}\, e^{-\lambda t_i})^{\theta}\}^{2}}{1-\{1-(\alpha e^{-\lambda t_i}/1-\overline{\alpha}\, e^{-\lambda t_i})^{\theta}\}^{a}}$$

$$\frac{\partial^2\ell}{\partial b^2}=-\frac{n}{b^2}$$

$$\frac{\partial^2 \ell}{\partial \lambda^2} = -\frac{n}{\lambda^2} + (\theta+1)\sum_{i=0}^{n} \frac{\bar{\alpha}^2 e^{-2\lambda t_i} t_i^2}{(1-\bar{\alpha} e^{-\lambda t_i})^2} + (\theta+1)\sum_{i=0}^{n} \frac{\bar{\alpha} e^{-\lambda t_i} t_i^2}{1-\bar{\alpha} e^{-\lambda t_i}}$$

$$+\theta^2 (1-a)\sum_{i=0}^{n} \frac{(\alpha e^{-\lambda t_i}/1-\bar{\alpha} e^{-\lambda t_i})^{2(\theta-1)} \left\{-\alpha\bar{\alpha} e^{-2\lambda t_i} t_i/(1-\bar{\alpha} e^{-\lambda t_i})^2 - (\alpha e^{-\lambda t_i} t_i/1-\bar{\alpha} e^{-\lambda t_i})\right\}^2}{\{1-(\alpha e^{-\lambda t_i}/1-\bar{\alpha} e^{-\lambda t_i})^{\theta}\}^2}$$

$$+\theta(\theta-1)(1-a)\sum_{i=0}^{n} \frac{(\alpha e^{-\lambda t_i}/1-\bar{\alpha} e^{-\lambda t_i})^{\theta-2} \left\{-\alpha\bar{\alpha} e^{-2\lambda t_i} t_i/(1-\bar{\alpha} e^{-\lambda t_i})^2 - (\alpha e^{-\lambda t_i} t_i/1-\bar{\alpha} e^{-\lambda t_i})\right\}^2}{1-(\alpha e^{-\lambda t_i}/1-\bar{\alpha} e^{-\lambda t_i})^{\theta}}$$

$$+\theta\,(1-a)\sum_{i=0}^{n} \frac{(\alpha e^{-\lambda t_i}/1-\bar{\alpha} e^{-\lambda t_i})^{\theta-1}}{1-(\alpha e^{-\lambda t_i}/1-\bar{\alpha} e^{-\lambda t_i})^{\theta}}$$

$$\times\left\{2\alpha\bar{\alpha}^2 e^{-3\lambda t_i} t_i^2/(1-\bar{\alpha} e^{-\lambda t_i})^3 + 3\alpha\bar{\alpha} e^{-2\lambda t_i} t_i^2/(1-\bar{\alpha} e^{-\lambda t_i})^2 + (\alpha e^{-\lambda t_i} t_i^2/1-\bar{\alpha} e^{-\lambda t_i})\right\}$$

$$+a^2\theta^2(1-b)\sum_{i=0}^{n} \frac{(\alpha e^{-\lambda t_i}/1-\bar{\alpha} e^{-\lambda t_i})^{2(\theta-1)} \{1-(\alpha e^{-\lambda t_i}/1-\bar{\alpha} e^{-\lambda t_i})^{\theta}\}^{2(a-1)}}{\{1-\{1-(\alpha e^{-\lambda t_i}/1-\bar{\alpha} e^{-\lambda t_i})^{\theta}\}^{a}\}^2}$$

$$\times\left\{-\alpha\bar{\alpha} e^{-2\lambda t_i} t_i/(1-\bar{\alpha} e^{-\lambda t_i})^2 - (\alpha e^{-\lambda t_i} t_i/1-\bar{\alpha} e^{-\lambda t_i})\right\}^2$$

$$+a(a-1)\theta^2(1-b)\sum_{i=0}^{n} \frac{(\alpha e^{-\lambda t_i}/1-\bar{\alpha} e^{-\lambda t_i})^{2(\theta-1)} \{1-(\alpha e^{-\lambda t_i}/1-\bar{\alpha} e^{-\lambda t_i})^{\theta}\}^{a-2}}{1-\{1-(\alpha e^{-\lambda t_i}/1-\bar{\alpha} e^{-\lambda t_i})^{\theta}\}^{a}}$$

$$\times\left\{-\alpha\bar{\alpha} e^{-2\lambda t_i} t_i/(1-\bar{\alpha} e^{-\lambda t_i})^2 - (\alpha e^{-\lambda t_i} t_i/1-\bar{\alpha} e^{-\lambda t_i})\right\}^2$$

$$+\theta\,(\theta-1)(1-b)\sum_{i=0}^{n} \frac{(\alpha e^{-\lambda t_i}/1-\bar{\alpha} e^{-\lambda t_i})^{\theta-2} \{1-(\alpha e^{-\lambda t_i}/1-\bar{\alpha} e^{-\lambda t_i})^{\theta}\}^{a-1}}{1-\{1-(\alpha e^{-\lambda t_i}/1-\bar{\alpha} e^{-\lambda t_i})^{\theta}\}^{a}}$$

$$\times\left\{-\alpha\bar{\alpha} e^{-2\lambda t_i} t_i/(1-\bar{\alpha} e^{-\lambda t_i})^2 - (\alpha e^{-\lambda t_i} t_i/1-\bar{\alpha} e^{-\lambda t_i})\right\}^2$$

$$+\theta\,(1-a)\sum_{i=0}^{n} \frac{(\alpha e^{-\lambda t_i}/1-\bar{\alpha} e^{-\lambda t_i})^{\theta-1} \{1-(\alpha e^{-\lambda t_i}/1-\bar{\alpha} e^{-\lambda t_i})^{\theta}\}^{a-1}}{1-(\alpha e^{-\lambda t_i}/1-\bar{\alpha} e^{-\lambda t_i})^{\theta}}$$

$$\times\left\{2\alpha\bar{\alpha}^2 e^{-3\lambda t_i} t_i^2/(1-\bar{\alpha} e^{-\lambda t_i})^3 + 3\alpha\bar{\alpha} e^{-2\lambda t_i} t_i^2/(1-\bar{\alpha} e^{-\lambda t_i})^2 + (\alpha e^{-\lambda t_i} t_i^2/1-\bar{\alpha} e^{-\lambda t_i})\right\}$$

$$\frac{\partial^2 \ell}{\partial \theta\,\partial \alpha} = \frac{n}{\alpha} - \sum_{i=0}^{n} \frac{e^{-\lambda t_i}}{1-\bar{\alpha} e^{-\lambda t_i}}$$

$$+(1-a)\sum_{i=0}^{n} \frac{(\alpha e^{-\lambda t_i}/1-\bar{\alpha} e^{-\lambda t_i})^{\theta-1} \left\{-\alpha e^{-2\lambda t_i}/(1-\bar{\alpha} e^{-\lambda t_i})^2 + (e^{-\lambda t_i}/1-\bar{\alpha} e^{-\lambda t_i})\right\}}{1-(\alpha e^{-\lambda t_i}/1-\bar{\alpha} e^{-\lambda t_i})^{\theta}}$$

$$+(1-a)\theta\sum_{i=0}^{n} \frac{(\alpha e^{-\lambda t_i}/1-\bar{\alpha} e^{-\lambda t_i})^{2\theta-1} \left\{-\alpha e^{-2\lambda t_i}/(1-\bar{\alpha} e^{-\lambda t_i})^2 + (e^{-\lambda t_i}/1-\bar{\alpha} e^{-\lambda t_i})\right\}\log\,(e^{-\lambda t_i}/1-\bar{\alpha} e^{-\lambda t_i})}{\{1-(\alpha e^{-\lambda t_i}/1-\bar{\alpha} e^{-\lambda t_i})^{\theta}\}^2}$$

$$+(1-a)\theta\sum_{i=0}^{n} \frac{(\alpha e^{-\lambda t_i}/1-\bar{\alpha} e^{-\lambda t_i})^{\theta-1} \left\{-\alpha e^{-2\lambda t_i}/(1-\bar{\alpha} e^{-\lambda t_i})^2 + (e^{-\lambda t_i}/1-\bar{\alpha} e^{-\lambda t_i})\right\}\log\,(e^{-\lambda t_i}/1-\bar{\alpha} e^{-\lambda t_i})}{1-(\alpha e^{-\lambda t_i}/1-\bar{\alpha} e^{-\lambda t_i})^{\theta}}$$

$$+a\,(b-1)\sum_{i=0}^{n} \frac{(\alpha e^{-\lambda t_i}/1-\bar{\alpha} e^{-\lambda t_i})^{\theta-1} \{1-(\alpha e^{-\lambda t_i}/1-\bar{\alpha} e^{-\lambda t_i})^{\theta}\}^{a-1}}{1-\{1-(\alpha e^{-\lambda t_i}/1-\bar{\alpha} e^{-\lambda t_i})^{\theta}\}^{a}}$$

$$\times\left\{-\alpha e^{-2\lambda t_i}/(1-\bar{\alpha} e^{-\lambda t_i})^2 + (e^{-\lambda t_i}/1-\bar{\alpha} e^{-\lambda t_i})\right\}$$

$$-a^{2}\theta(b-1)\sum_{i=0}^{n}\frac{(\alpha e^{-\lambda t_i}/1-\bar{\alpha} e^{-\lambda t_i})^{2\theta-1}\{1-(\alpha e^{-\lambda t_i}/1-\bar{\alpha} e^{-\lambda t_i})^{\theta}\}^{2(a-1)}}{[1-\{1-(\alpha e^{-\lambda t_i}/1-\bar{\alpha} e^{-\lambda t_i})^{\theta}\}^{a}]^{2}}$$

$$\times\left\{-\alpha e^{-2\lambda t_i}/(1-\bar{\alpha} e^{-\lambda t_i})^{2}+(e^{-\lambda t_i}/1-\bar{\alpha} e^{-\lambda t_i})\right\}\log(e^{-\lambda t_i}/1-\bar{\alpha} e^{-\lambda t_i})$$

$$-a(a-1)\theta(b-1)\sum_{i=0}^{n}\frac{(\alpha e^{-\lambda t_i}/1-\bar{\alpha} e^{-\lambda t_i})^{2\theta-1}\{1-(\alpha e^{-\lambda t_i}/1-\bar{\alpha} e^{-\lambda t_i})^{\theta}\}^{a-2}}{1-\{1-(\alpha e^{-\lambda t_i}/1-\bar{\alpha} e^{-\lambda t_i})^{\theta}\}^{a}}$$

$$\times\left\{-\alpha e^{-2\lambda t_i}/(1-\bar{\alpha} e^{-\lambda t_i})^{2}+(e^{-\lambda t_i}/1-\bar{\alpha} e^{-\lambda t_i})\right\}\log(e^{-\lambda t_i}/1-\bar{\alpha} e^{-\lambda t_i})$$

$$+a\,\theta(b-1)\sum_{i=0}^{n}\frac{(\alpha e^{-\lambda t_i}/1-\bar{\alpha} e^{-\lambda t_i})^{\theta-1}\{1-(\alpha e^{-\lambda t_i}/1-\bar{\alpha} e^{-\lambda t_i})^{\theta}\}^{a-1}}{1-\{1-(\alpha e^{-\lambda t_i}/1-\bar{\alpha} e^{-\lambda t_i})^{\theta}\}^{a}}$$

$$\times\left\{-\alpha e^{-2\lambda t_i}/(1-\bar{\alpha} e^{-\lambda t_i})^{2}+(e^{-\lambda t_i}/1-\bar{\alpha} e^{-\lambda t_i})\right\}\log(e^{-\lambda t_i}/1-\bar{\alpha} e^{-\lambda t_i})$$

$$\frac{\partial^{2}\ell}{\partial\theta\,\partial a}=-\sum_{i=0}^{n}\frac{(\alpha e^{-\lambda t_i}/1-\bar{\alpha} e^{-\lambda t_i})^{\theta}\log(\alpha e^{-\lambda t_i}/1-\bar{\alpha} e^{-\lambda t_i})}{1-(\alpha e^{-\lambda t_i}/1-\bar{\alpha} e^{-\lambda t_i})^{\theta}}$$

$$+(b-1)\sum_{i=0}^{n}\frac{(\alpha e^{-\lambda t_i}/1-\bar{\alpha} e^{-\lambda t_i})^{\theta}\,[1-(\alpha e^{-\lambda t_i}/1-\bar{\alpha} e^{-\lambda t_i})^{\theta}]^{a-1}\log(\alpha e^{-\lambda t_i}/1-\bar{\alpha} e^{-\lambda t_i})}{1-[1-(\alpha e^{-\lambda t_i}/1-\bar{\alpha} e^{-\lambda t_i})^{\theta}]^{a}}$$

$$+a(b-1)\sum_{i=0}^{n}\frac{(\alpha e^{-\lambda t_i}/1-\bar{\alpha} e^{-\lambda t_i})^{\theta}\,[1-(\alpha e^{-\lambda t_i}/1-\bar{\alpha} e^{-\lambda t_i})^{\theta}]^{2a-1}\log(\alpha e^{-\lambda t_i}/1-\bar{\alpha} e^{-\lambda t_i})\log[1-(\alpha e^{-\lambda t_i}/1-\bar{\alpha} e^{-\lambda t_i})^{\theta}]}{[1-[1-(\alpha e^{-\lambda t_i}/1-\bar{\alpha} e^{-\lambda t_i})^{\theta}]^{a}]^{2}}$$

$$+a(b-1)\sum_{i=0}^{n}\frac{(\alpha e^{-\lambda t_i}/1-\bar{\alpha} e^{-\lambda t_i})^{\theta}\,[1-(\alpha e^{-\lambda t_i}/1-\bar{\alpha} e^{-\lambda t_i})^{\theta}]^{a-1}\log(\alpha e^{-\lambda t_i}/1-\bar{\alpha} e^{-\lambda t_i})\log[1-(\alpha e^{-\lambda t_i}/1-\bar{\alpha} e^{-\lambda t_i})^{\theta}]}{1-[1-(\alpha e^{-\lambda t_i}/1-\bar{\alpha} e^{-\lambda t_i})^{\theta}]^{a}}$$

$$\frac{\partial^{2}\ell}{\partial\theta\,\partial b}=\sum_{i=0}^{n}\frac{a[1-(\alpha e^{-\lambda t_i}/1-\bar{\alpha} e^{-\lambda t_i})^{\theta}]^{a-1}(\alpha e^{-\lambda t_i}/1-\bar{\alpha} e^{-\lambda t_i})^{\theta}\log(\alpha e^{-\lambda t_i}/1-\bar{\alpha} e^{-\lambda t_i})}{1-[1-(\alpha e^{-\lambda t_i}/1-\bar{\alpha} e^{-\lambda t_i})^{\theta}]^{a}}$$

$$\frac{\partial^{2}\ell}{\partial\theta\,\partial\lambda}=-\sum_{i=0}^{n}t_i-\sum_{i=0}^{n}\frac{\bar{\alpha} e^{-\lambda t_i}\,t_i}{1-\bar{\alpha} e^{-\lambda t_i}}$$

$$+(1-a)\sum_{i=0}^{n}\frac{(\alpha e^{-\lambda t_i}/1-\bar{\alpha} e^{-\lambda t_i})^{\theta-1}\left\{-\alpha\bar{\alpha} e^{-2\lambda t_i}\,t_i/(1-\bar{\alpha} e^{-\lambda t_i})^{2}-(\alpha e^{-\lambda t_i}\,t_i/1-\bar{\alpha} e^{-\lambda t_i})\right\}}{1-(\alpha e^{-\lambda t_i}/1-\bar{\alpha} e^{-\lambda t_i})^{\theta}}$$

$$+(1-a)\theta\sum_{i=0}^{n}\frac{(\alpha e^{-\lambda t_i}/1-\bar{\alpha} e^{-\lambda t_i})^{2\theta-1}\left\{-\alpha\bar{\alpha} e^{-2\lambda t_i}\,t_i/(1-\bar{\alpha} e^{-\lambda t_i})^{2}-(\alpha e^{-\lambda t_i}\,t_i/1-\bar{\alpha} e^{-\lambda t_i})\right\}}{\{1-(\alpha e^{-\lambda t_i}/1-\bar{\alpha} e^{-\lambda t_i})^{\theta}\}^{2}}$$

$$\times\log(\alpha e^{-\lambda t_i}\,t_i/1-\bar{\alpha} e^{-\lambda t_i})$$

$$+(1-a)\theta\sum_{i=0}^{n}\frac{(\alpha e^{-\lambda t_i}/1-\bar{\alpha} e^{-\lambda t_i})^{\theta-1}\left\{-\alpha\bar{\alpha} e^{-2\lambda t_i}\,t_i/(1-\bar{\alpha} e^{-\lambda t_i})^{2}-(\alpha e^{-\lambda t_i}\,t_i/1-\bar{\alpha} e^{-\lambda t_i})\right\}}{1-(\alpha e^{-\lambda t_i}/1-\bar{\alpha} e^{-\lambda t_i})^{\theta}}$$

$$\times\log(\alpha e^{-\lambda t_i}\,t_i/1-\bar{\alpha} e^{-\lambda t_i})$$

$$+a\,(b-1)\sum_{i=0}^{n}\frac{(\alpha e^{-\lambda t_i}/1-\bar{\alpha} e^{-\lambda t_i})^{\theta-1}\{1-(\alpha e^{-\lambda t_i}/1-\bar{\alpha} e^{-\lambda t_i})^{\theta}\}^{a-1}}{1-\{1-(\alpha e^{-\lambda t_i}/1-\bar{\alpha} e^{-\lambda t_i})^{\theta}\}^{a}}$$

$$\times\left\{-\alpha\bar{\alpha} e^{-2\lambda t_i}\,t_i/(1-\bar{\alpha} e^{-\lambda t_i})^{2}-(\alpha e^{-\lambda t_i}\,t_i/1-\bar{\alpha} e^{-\lambda t_i})\right\}$$

$$+ a^2\theta(b-1)\sum_{i=0}^{n}\frac{(\alpha e^{-\lambda t_i}/1-\bar{\alpha} e^{-\lambda t_i})^{2\theta-1}\{1-(\alpha e^{-\lambda t_i}/1-\bar{\alpha} e^{-\lambda t_i})^{\theta}\}^{2(a-1)}}{[1-\{1-(\alpha e^{-\lambda t_i}/1-\bar{\alpha} e^{-\lambda t_i})^{\theta}\}^{a}]^2}$$

$$\times\left\{-\alpha\bar{\alpha} e^{-2\lambda t_i} t_i/(1-\bar{\alpha} e^{-\lambda t_i})^2 - (\alpha e^{-\lambda t_i} t_i/1-\bar{\alpha} e^{-\lambda t_i})\right\}\log(\alpha e^{-\lambda t_i} t_i/1-\bar{\alpha} e^{-\lambda t_i})$$

$$+ a(a-1)\theta(b-1)\sum_{i=0}^{n}\frac{(\alpha e^{-\lambda t_i}/1-\bar{\alpha} e^{-\lambda t_i})^{2\theta-1}\{1-(\alpha e^{-\lambda t_i}/1-\bar{\alpha} e^{-\lambda t_i})^{\theta}\}^{a-2}}{1-\{1-(\alpha e^{-\lambda t_i}/1-\bar{\alpha} e^{-\lambda t_i})^{\theta}\}^{a}}$$

$$\times\left\{-\alpha\bar{\alpha} e^{-2\lambda t_i} t_i/(1-\bar{\alpha} e^{-\lambda t_i})^2 - (\alpha e^{-\lambda t_i} t_i/1-\bar{\alpha} e^{-\lambda t_i})\right\}\log(\alpha e^{-\lambda t_i} t_i/1-\bar{\alpha} e^{-\lambda t_i})$$

$$+ a\,\theta(b-1)\sum_{i=0}^{n}\frac{(\alpha e^{-\lambda t_i}/1-\bar{\alpha} e^{-\lambda t_i})^{\theta-1}\{1-(\alpha e^{-\lambda t_i}/1-\bar{\alpha} e^{-\lambda t_i})^{\theta}\}^{a-1}}{1-\{1-(\alpha e^{-\lambda t_i}/1-\bar{\alpha} e^{-\lambda t_i})^{\theta}\}^{a}}$$

$$\times\left\{-\alpha\bar{\alpha} e^{-2\lambda t_i} t_i/(1-\bar{\alpha} e^{-\lambda t_i})^2 - (\alpha e^{-\lambda t_i} t_i/1-\bar{\alpha} e^{-\lambda t_i})\right\}\log(\alpha e^{-\lambda t_i} t_i/1-\bar{\alpha} e^{-\lambda t_i})$$

$$\frac{\partial^2\ell}{\partial\alpha\,\partial a} = \theta\sum_{i=0}^{n} - \frac{(\alpha e^{-\lambda t_i}/1-\bar{\alpha} e^{-\lambda t_i})^{\theta-1}\left\{-\alpha\bar{\alpha} e^{-2\lambda t_i} t_i/(1-\bar{\alpha} e^{-\lambda t_i})^2 + (\alpha e^{-\lambda t_i} t_i/1-\bar{\alpha} e^{-\lambda t_i})\right\}}{1-(\alpha e^{-\lambda t_i}/1-\bar{\alpha} e^{-\lambda t_i})^{\theta}}$$

$$+ (b-1)\,\theta\sum_{i=0}^{n}\frac{(\alpha e^{-\lambda t_i}/1-\bar{\alpha} e^{-\lambda t_i})^{\theta-1}\{1-(\alpha e^{-\lambda t_i}/1-\bar{\alpha} e^{-\lambda t_i})^{\theta}\}^{a-1}}{1-\{1-(\alpha e^{-\lambda t_i}/1-\bar{\alpha} e^{-\lambda t_i})^{\theta}\}^{a}}$$

$$\times\left\{-\alpha\bar{\alpha} e^{-2\lambda t_i} t_i/(1-\bar{\alpha} e^{-\lambda t_i})^2 + (\alpha e^{-\lambda t_i} t_i/1-\bar{\alpha} e^{-\lambda t_i})\right\}$$

$$+ (b-1)\,a\theta\sum_{i=0}^{n}\frac{(\alpha e^{-\lambda t_i}/1-\bar{\alpha} e^{-\lambda t_i})^{\theta-1}\{1-(\alpha e^{-\lambda t_i}/1-\bar{\alpha} e^{-\lambda t_i})^{\theta}\}^{2a-1}}{[1-\{1-(\alpha e^{-\lambda t_i}/1-\bar{\alpha} e^{-\lambda t_i})^{\theta}\}^{a}]^2}$$

$$\times\left\{-\alpha\bar{\alpha} e^{-2\lambda t_i} t_i/(1-\bar{\alpha} e^{-\lambda t_i})^2 + (\alpha e^{-\lambda t_i} t_i/1-\bar{\alpha} e^{-\lambda t_i})\right\}\log[1-(\alpha e^{-\lambda t_i}/1-\bar{\alpha} e^{-\lambda t_i})^{\theta}]$$

$$+ (b-1)\,a\theta\sum_{i=0}^{n}\frac{(\alpha e^{-\lambda t_i}/1-\bar{\alpha} e^{-\lambda t_i})^{\theta-1}\{1-(\alpha e^{-\lambda t_i}/1-\bar{\alpha} e^{-\lambda t_i})^{\theta}\}^{a-1}}{1-\{1-(\alpha e^{-\lambda t_i}/1-\bar{\alpha} e^{-\lambda t_i})^{\theta}\}^{a}}$$

$$\times\left\{-\alpha\bar{\alpha} e^{-2\lambda t_i} t_i/(1-\bar{\alpha} e^{-\lambda t_i})^2 + (\alpha e^{-\lambda t_i} t_i/1-\bar{\alpha} e^{-\lambda t_i})\right\}\log[1-(\alpha e^{-\lambda t_i}/1-\bar{\alpha} e^{-\lambda t_i})^{\theta}]$$

$$\frac{\partial^2\ell}{\partial\alpha\,\partial b} = a\theta\sum_{i=0}^{n}\frac{(\alpha e^{-\lambda t_i}/1-\bar{\alpha} e^{-\lambda t_i})^{\theta-1}\{1-(\alpha e^{-\lambda t_i}/1-\bar{\alpha} e^{-\lambda t_i})^{\theta}\}^{a-1}}{1-\{1-(\alpha e^{-\lambda t_i}/1-\bar{\alpha} e^{-\lambda t_i})^{\theta}\}^{a}}$$

$$\times\left\{-\alpha\bar{\alpha} e^{-2\lambda t_i} t_i/(1-\bar{\alpha} e^{-\lambda t_i})^2 + (\alpha e^{-\lambda t_i} t_i/1-\bar{\alpha} e^{-\lambda t_i})\right\}$$

$$\frac{\partial^2\ell}{\partial\alpha\,\partial\lambda} = (\theta+1)\sum_{i=0}^{n}\left\{\bar{\alpha} e^{-2\lambda t_i} t_i/(1-\bar{\alpha} e^{-\lambda t_i})^2 + (e^{-\lambda t_i} t_i/1-\bar{\alpha} e^{-\lambda t_i})\right\}$$

$$+(1-a)\theta\sum_{i=0}^{n}\frac{(\alpha e^{-\lambda t_i}/1-\bar{\alpha} e^{-\lambda t_i})^{\theta-1}}{1-(\alpha e^{-\lambda t_i}/1-\bar{\alpha} e^{-\lambda t_i})^{\theta}}$$

$$\times\left\{2\alpha\bar{\alpha} e^{-3\lambda t_i} t_i/(1-\bar{\alpha} e^{-\lambda t_i})^3 - \{\bar{\alpha} e^{-2\lambda t_i} t_i/(1-\bar{\alpha} e^{-\lambda t_i})^2\} + 2\alpha\, e^{-2\lambda t_i} t_i/(1-\bar{\alpha} e^{-\lambda t_i})^2 - (\alpha e^{-\lambda t_i} t_i/1-\bar{\alpha} e^{-\lambda t_i})\right\}$$

$$+(1-a)\;\theta^2\sum_{i=0}^{n}\frac{(\alpha e^{-\lambda t_i}/1-\bar{\alpha} e^{-\lambda t_i})^{2(\theta-1)}}{\{1-(\alpha e^{-\lambda t_i}/1-\bar{\alpha} e^{-\lambda t_i})^{\theta}\}^2}$$

$$\times[-\alpha e^{-2\lambda t_i} t_i/(1-\bar{\alpha} e^{-\lambda t_i})^2 + \{e^{-2\lambda t_i}/(1-\bar{\alpha} e^{-\lambda t_i})\}]\times[-\alpha\bar{\alpha} e^{-2\lambda t_i} t_i/(1-\bar{\alpha} e^{-\lambda t_i})^2 - (\alpha e^{-\lambda t_i} t_i/1-\bar{\alpha} e^{-\lambda t_i})]$$

$$+\ \theta(\theta-1)\sum_{i=0}^{n}\frac{(\alpha e^{-\lambda t_i}/1-\bar{\alpha}e^{-\lambda t_i})^{\theta-2}}{1-(\alpha e^{-\lambda t_i}/1-\bar{\alpha}e^{-\lambda t_i})^{\theta}}$$

$$\times[-\alpha e^{-2\lambda t_i}t_i/(1-\bar{\alpha}e^{-\lambda t_i})^2+\{e^{-2\lambda t_i}/(1-\bar{\alpha}e^{-\lambda t_i})\}]\times[-\alpha\bar{\alpha}e^{-2\lambda t_i}t_i/(1-\bar{\alpha}e^{-\lambda t_i})^2-(\alpha e^{-\lambda t_i}t_i/1-\bar{\alpha}e^{-\lambda t_i})]$$

$$+a\ \theta\sum_{i=0}^{n}\frac{(\alpha e^{-\lambda t_i}/1-\bar{\alpha}e^{-\lambda t_i})^{\theta-1}\{1-(\alpha e^{-\lambda t_i}/1-\bar{\alpha}e^{-\lambda t_i})^{\theta}\}^{a-1}}{1-\{1-(\alpha e^{-\lambda t_i}/1-\bar{\alpha}e^{-\lambda t_i})^{\theta}\}^{a}}$$

$$\times\left\{2\alpha\bar{\alpha}e^{-3\lambda t_i}t_i/(1-\bar{\alpha}e^{-\lambda t_i})^3-\{\bar{\alpha}e^{-2\lambda t_i}t_i/(1-\bar{\alpha}e^{-\lambda t_i})^2\}+2\alpha\ e^{-2\lambda t_i}t_i/(1-\bar{\alpha}e^{-\lambda t_i})^2-(\alpha e^{-\lambda t_i}t_i/1-\bar{\alpha}e^{-\lambda t_i})\right\}$$

$$-a^2\ \theta^2\sum_{i=0}^{n}\frac{(\alpha e^{-\lambda t_i}/1-\bar{\alpha}e^{-\lambda t_i})^{2(\theta-1)}\{1-(\alpha e^{-\lambda t_i}/1-\bar{\alpha}e^{-\lambda t_i})^{\theta}\}^{2(a-1)}}{[1-\{1-(\alpha e^{-\lambda t_i}/1-\bar{\alpha}e^{-\lambda t_i})^{\theta}\}^{a}]^2}$$

$$\times[-\alpha e^{-2\lambda t_i}t_i/(1-\bar{\alpha}e^{-\lambda t_i})^2+\{e^{-2\lambda t_i}/(1-\bar{\alpha}e^{-\lambda t_i})\}]\times[-\alpha\bar{\alpha}e^{-2\lambda t_i}t_i/(1-\bar{\alpha}e^{-\lambda t_i})^2-(\alpha e^{-\lambda t_i}t_i/1-\bar{\alpha}e^{-\lambda t_i})]$$

$$-a(a-1)\ \theta^2\sum_{i=0}^{n}\frac{(\alpha e^{-\lambda t_i}/1-\bar{\alpha}e^{-\lambda t_i})^{2(\theta-1)}\{1-(\alpha e^{-\lambda t_i}/1-\bar{\alpha}e^{-\lambda t_i})^{\theta}\}^{a-2}}{1-\{1-(\alpha e^{-\lambda t_i}/1-\bar{\alpha}e^{-\lambda t_i})^{\theta}\}^{a}}$$

$$\times[-\alpha e^{-2\lambda t_i}t_i/(1-\bar{\alpha}e^{-\lambda t_i})^2+\{e^{-2\lambda t_i}/(1-\bar{\alpha}e^{-\lambda t_i})\}]\times[-\alpha\bar{\alpha}e^{-2\lambda t_i}t_i/(1-\bar{\alpha}e^{-\lambda t_i})^2-(\alpha e^{-\lambda t_i}t_i/1-\bar{\alpha}e^{-\lambda t_i})]$$

$$+a\ \theta(\theta-1)\sum_{i=0}^{n}\frac{(\alpha e^{-\lambda t_i}/1-\bar{\alpha}e^{-\lambda t_i})^{\theta-2}\{1-(\alpha e^{-\lambda t_i}/1-\bar{\alpha}e^{-\lambda t_i})^{\theta}\}^{a-1}}{1-\{1-(\alpha e^{-\lambda t_i}/1-\bar{\alpha}e^{-\lambda t_i})^{\theta}\}^{a}}$$

$$\times[-\alpha e^{-2\lambda t_i}t_i/(1-\bar{\alpha}e^{-\lambda t_i})^2+\{e^{-2\lambda t_i}/(1-\bar{\alpha}e^{-\lambda t_i})\}]\times[-\alpha\bar{\alpha}e^{-2\lambda t_i}t_i/(1-\bar{\alpha}e^{-\lambda t_i})^2-(\alpha e^{-\lambda t_i}t_i/1-\bar{\alpha}e^{-\lambda t_i})]$$

$$\frac{\partial^2\ell}{\partial a\,\partial b}=-\sum_{i=0}^{n}\frac{\{1-(\alpha e^{-\lambda t_i}/1-\bar{\alpha}e^{-\lambda t_i})^{\theta}\}^{a}\log[1-(\alpha e^{-\lambda t_i}/1-\bar{\alpha}e^{-\lambda t_i})^{\theta}]}{1-\{1-(\alpha e^{-\lambda t_i}/1-\bar{\alpha}e^{-\lambda t_i})^{\theta}\}^{a}}$$

$$\frac{\partial^2\ell}{\partial a\,\partial\lambda}=-\theta\sum_{i=0}^{n}\frac{(\alpha e^{-\lambda t_i}/1-\bar{\alpha}e^{-\lambda t_i})^{\theta-1}[-\alpha\bar{\alpha}e^{-2\lambda t_i}t_i/(1-\bar{\alpha}e^{-\lambda t_i})^2-(\alpha e^{-\lambda t_i}t_i/1-\bar{\alpha}e^{-\lambda t_i})]}{1-(\alpha e^{-\lambda t_i}/1-\bar{\alpha}e^{-\lambda t_i})^{\theta}}$$

$$+(b-1)\ \theta\sum_{i=0}^{n}\frac{(\alpha e^{-\lambda t_i}/1-\bar{\alpha}e^{-\lambda t_i})^{\theta-1}\{1-(\alpha e^{-\lambda t_i}/1-\bar{\alpha}e^{-\lambda t_i})^{\theta}\}^{a-1}}{1-\{1-(\alpha e^{-\lambda t_i}/1-\bar{\alpha}e^{-\lambda t_i})^{\theta}\}^{a}}$$

$$\times[-\alpha\bar{\alpha}e^{-2\lambda t_i}t_i/(1-\bar{\alpha}e^{-\lambda t_i})^2-(\alpha e^{-\lambda t_i}t_i/1-\bar{\alpha}e^{-\lambda t_i})]$$

$$+(b-1)a\theta\sum_{i=0}^{n}\frac{(\alpha e^{-\lambda t_i}/1-\bar{\alpha}e^{-\lambda t_i})^{\theta-1}\{1-(\alpha e^{-\lambda t_i}/1-\bar{\alpha}e^{-\lambda t_i})^{\theta}\}^{2a-1}}{[1-\{1-(\alpha e^{-\lambda t_i}/1-\bar{\alpha}e^{-\lambda t_i})^{\theta}\}^{a}]^2}$$

$$\times[-\alpha\bar{\alpha}e^{-2\lambda t_i}t_i/(1-\bar{\alpha}e^{-\lambda t_i})^2-(\alpha e^{-\lambda t_i}t_i/1-\bar{\alpha}e^{-\lambda t_i})]\log[1-(\alpha e^{-\lambda t_i}/1-\bar{\alpha}e^{-\lambda t_i})^{\theta}]$$

$$+(b-1)a\theta\sum_{i=0}^{n}\frac{(\alpha e^{-\lambda t_i}/1-\bar{\alpha}e^{-\lambda t_i})^{\theta-1}\{1-(\alpha e^{-\lambda t_i}/1-\bar{\alpha}e^{-\lambda t_i})^{\theta}\}^{a-1}}{1-\{1-(\alpha e^{-\lambda t_i}/1-\bar{\alpha}e^{-\lambda t_i})^{\theta}\}^{a}}$$

$$\times[-\alpha\bar{\alpha}e^{-2\lambda t_i}t_i/(1-\bar{\alpha}e^{-\lambda t_i})^2-(\alpha e^{-\lambda t_i}t_i/1-\bar{\alpha}e^{-\lambda t_i})]\log[1-(\alpha e^{-\lambda t_i}/1-\bar{\alpha}e^{-\lambda t_i})^{\theta}]$$

$$\frac{\partial^2\ell}{\partial b\,\partial\lambda}=a\theta\sum_{i=0}^{n}\frac{(\alpha e^{-\lambda t_i}/1-\bar{\alpha}e^{-\lambda t_i})^{\theta-1}\{1-(\alpha e^{-\lambda t_i}/1-\bar{\alpha}e^{-\lambda t_i})^{\theta}\}^{a-1}}{1-\{1-(\alpha e^{-\lambda t_i}/1-\bar{\alpha}e^{-\lambda t_i})^{\theta}\}^{a}}$$

$$\times[-\alpha\bar{\alpha}e^{-2\lambda t_i}t_i/(1-\bar{\alpha}e^{-\lambda t_i})^2-(\alpha e^{-\lambda t_i}t_i/1-\bar{\alpha}e^{-\lambda t_i})]$$

Where $\psi'(.)$ is the derivative of the digamma function.